\documentstyle[11pt,leqno,amscd,amssymb,pstricks,verbatim]{amsart}
\newtheorem{thm}{Theorem}[section]
\newtheorem{lem}[thm]{Lemma}
\newtheorem{cor}[thm]{Corollary}
\newtheorem{prop}[thm]{Proposition}

\theoremstyle{definition}
\newtheorem{example}[thm]{Example}

\newtheorem{defn}[thm]{Definition}

\newtheorem{rem}[thm]{Remark}

\numberwithin{equation}{thm}
\def\lan{\langle}    \def\ran{\rangle}
\def\lr#1{\langle #1\rangle}

\def\az{\alpha}  
\def\bz{\beta}  \def\dt{\Delta}
\def\gz{\gamma}
  \def\ooz{\Omega}
\def\sz{\sigma}

\def\vez{\varepsilon}  
 \def\dz{\delta}

\def\lz{\lambda}
\def\cm{{\cal M}}\def\cp{{\cal P}} \def\fr{{\frak R}}      \def\cc{{\cal C}} \def\cx{{\cal X}} \def\cy{{\cal Y}}
\def\ch{{\cal H}}     \def\fs{{\frak S}}
\def\cz{{\cal Z}} \def\co{{\cal O}} \def\ce{{\cal E}}
  
\def\ca{{\cal A}}
\def\bbn{{\mathbb N}}  \def\bbz{{\mathbb Z}}  \def\bbq{{\mathbb Q}}
  \def\bbt{{\mathbb T}}  
     \def\tpi{{\tilde \pi}}    \def\tlz{{\tilde \lz}}
 \def\vphi{\varphi}
\def\dpeq{{\prec\!\!\!\cdot}}        \def\bd{{\bf d}} \def\bU{{\bf U}}  
\def\peq{\preccurlyeq}  \def\leq{\leqslant}  \def\geq{\geqslant}
\def\ged{\geqslant} \def\led{\leqslant}
\def\lra{\longrightarrow}   \def\skrel{\stackrel}
\def\lmto{\longmapsto}
\def\ra{\rightarrow}
\def\hom{\mbox{\rm Hom}}

\def\ext{\mbox{\rm Ext}\,}   \def\Ker{\mbox{\rm Ker}\,}
\def\dim{\mbox{\rm dim}\,}
\def\udim{{\mathbf dim}} 
  
\def\rad{\mbox{\rm rad}\,}  \def\top{\mbox{\rm top}\,}
\def\rep{\mbox{\rm Rep}\,}

\begin{document}
\title[Monomial Bases for quantum affine
${\frak sl}_n$]%
{Monomial Bases for quantum affine
${\frak sl}_n$}
\author{Bangming Deng and Jie Du}
\address{Department of Mathematics, Beijing Normal University,
Beijing 100875,  China.} \email{dengbm@@bnu.edu.cn}
\address{School of Mathematics, University of New South Wales,
Sydney 2052, Australia.} \email{j.du@@unsw.edu.au\quad {\it Home
page}: http://www.maths.unsw.edu.au/$\sim$jied}

\thanks{Supported partially by ARC. The research was carried out
while the first author was visiting the University of New South Wales.
The hospitality and support of UNSW are gratefully acknowledged. The
first author is partially supported by the National Natural Science
Foundation of China.}


\subjclass[2000]{17B37, 16G20}

\begin{abstract}
We use the idea of generic extensions to investigate the correspondence
between the isomorphism classes of
nilpotent representations of a cyclic quiver and the orbits in the
corresponding representation varieties. We endow
the set $\cm$ of such isoclasses with a monoid structure and identify
the submonoid $\cm_c$ generated by simple modules.
On the other hand, we use the partial
ordering on the orbits (i.e., the Bruhat-Chevalley type ordering) to
induce a poset structure on $\cm$ and describe
the poset ideals generated by an element of the submonoid $\cm_c$
 in terms of the existence of a certain composition
series of the corresponding module.
As applications of these results,
we generalize some results of Ringel involving
special words to results with no restriction on words and
obtain a systematic
description of many monomial bases for any given
quantum affine ${\frak {sl}}_n$.
\end{abstract}

\maketitle

\section{Introduction}

Let $\bU$ be a quantum group over $\bbq(v)$ associated to a Cartan
datum $(I,\centerdot)$ in the sense of \cite[1.1]{L93}, and let
$E_i, F_i, K_i^{\pm1}$ ($i\in I$) be its generators. Then all
monomials in the $E_i$'s, $F_i$'s and $K_i^{\pm 1}$'s span $\bU$.
It is natural to ask which monomials form bases for $\bU$. Since
$\bU$ admits a triangular decomposition
$$\bU=\bU^-\bU^0\bU^+$$
where $\bU^+$ (resp. $\bU^-$, $\bU^0$) is the subalgebra generated
by the $E_i$'s (resp. $F_i$'s, $K_i^{\pm 1}$'s), and the monomials
$K_{\mathbf a}=K_1^{a_1}\cdots K_n^{a_n}$ with ${\mathbf a}
=\sum_ia_ii\in\bbz I$ form a basis for $\bU^0$, it would be
interesting to know monomial bases for $\bU^+$ (and hence for
$\bU^-$). More precisely, let $\ooz$ be the set of words on the
alphabet $I$. For each word $w=i_1\cdots i_m\in\ooz$, let
$$E_w=E_{i_1}\cdots E_{i_m},\quad F_w=F_{i_1}\cdots F_{i_m}.$$
We are interested in finding subsets $\ooz'\subset\ooz$ such that
the set $\{E_w\}_{w\in\ooz'}$ forms a basis for $\bU^+$.

In the case where $(I,\centerdot)$ is of simply-laced finite type,
Lusztig introduced certain monomial bases for $\bU^+$ in
\cite[7.8]{L90} relative to a fixed reduced expression of the
longest word in the corresponding Weyl group. Interesting
applications of monomial bases can be found in, e.g.,  \cite{R95},
\cite{CX}, \cite{Re1} and \cite{DP}. In \cite{DP} it is proved
that some (integral) monomial bases for a quantum $\frak{gl}_n$
give rise to monomial bases for $q$-Schur algebras and Hecke
algebras (see, e.g., \cite[7.2,9.4]{DP}).

In the investigation \cite{R93} of realizing the generic composition
algebra of a cyclic quiver as the $+$part of a quantum affine $\frak {sl}_n$,
Ringel constructed certain monomial bases over some so-called
condensed words
whose definition is rather long and complicated.
In this paper, we shall describe monomial bases for a quantum
affine ${\frak sl}_n$ in a more general and
satisfactory way. We shall prove the following monomial basis theorem
(see \S8).

\begin{thm}\label{MBT} Let $\bU=\bU_v(\widehat{\frak {sl}}_n)$. Then there is a partition
$\ooz=\cup_{\pi\in\Pi^s}\ooz_\pi$ such that, for any subsets $\ooz^\pm\subset\ooz$
with $|\ooz^\pm\cap \ooz_\pi|=1$ for every $\pi\in\Pi^s$, the set
$$\{F_yK_{\mathbf a}E_w\mid y\in\ooz^-,w\in\ooz^+,{\mathbf a}\in\Bbb Z I\}$$
forms a basis for $\bU$.
\end{thm}

The proof of the theorem uses the idea of generic extensions
of nilpotent representations of a cyclic quiver with $n$ vertices. We first
generalize a recent work by Reineke \cite{Re} to obtain a monoid
structure on the set $\cm$
 of isoclasses of nilpotent representations indexed by $\Pi$, the set of
$n$-tuples of partitions.
Simple modules will generate a  submonoid $\cm_c$. Thus, every word
$w\in\ooz$ defines a unique isoclass in $\cm_c$ which is indexed by $\wp(w)$.
We shall prove that the set of all $\wp(w)$ coincides with
the subset $\Pi^s$ of separated multipartitions defined in \cite[4.1]{R93}.
Thus, the fibres of $\wp$ yield a partition of $\ooz$.
Our argument will then depend heavily on the structure of nilpotent representations.

We organise the paper as follows. We start with investigating
several useful properties of nilpotent representations in \S2.
Then we move on looking at generic extensions of nilpotent
representations and constructing the monoid $\cm$ in \S3. In \S 4,
we discuss the relation between the submonoid $\cm_c$ and
separated multipartitions. Some algorithms are introduced to
calculate the multipartition $\wp(w)$ and the words in the fibres
of $\wp$. The notion of distinguished words is defined  in terms
of a certain module structure. However a combinatorial criterion
exists. These will be discussed in \S 5. There are two partial
order relations on nilpotent representations defined geometrically
by the inclusive relation on the closures of orbits and
algebraically by module extensions. We shall prove in \S6 that the
two relations coincide. Moreover, we describe a poset ideal
generated by an element of $\cm_c$ in terms of the existence of
certain composition series. This is Theorem \ref{DOMI}. Some
applications are given in the last three sections. In \S7, we
generalize some results of Ringel given in \cite{R93}, and
in \S 8, we prove Theorem \ref{MBT}.
Finally, in the last section, we construct a PBW type basis for
$\bU^+$ and speculate some relations between the various bases
including the canonical basis.

In a forthcoming paper, we expect to prove a similar result for the
quantum groups of simply-laced finite type.

\noindent
{\bf Some notation.}
Throughout this paper, $k$ denotes a field. We shall assume
that $k=\bar k$ is algebraically closed in section three.
All modules $M$ are finite dimensional over $k$. We denote by $\rad (M)$ the radical
of $M$, i.e., the intersection of all maximal submodules of $M$, and by
$\top (M):=M/\rad(M)$ the top of $M$.

The following lemma is a result of a pull-back or push-out diagram.

\begin{lem} \label{pbk}
Let $0\to X\to M\to Y\to 0$ be a short exact sequence of modules.
Then any exact sequence $0\to Y'\to Y\to Y''\to 0$ gives rise to
a commutative diagram with exact rows and columns

{\small
\begin{center}
\begin{pspicture}(0,.2)(7,4.1)
\psset{xunit=.7cm,yunit=.6cm,linewidth=0.5pt}
\uput[u](4,6){0}
\uput[u](6,6){0}
\psline{->}(4,6)(4,5.4)
\psline{->}(6,6)(6,5.4)

\uput[l](0.5,5){0}
\uput[l](2.5,5){$X$}
\uput[l](4.5,5){$M'$}
\uput[l](6.5,5){$Y'$}
\uput[l](8.15,5){$0$}
\psline{->}(0.6,5)(1.3,5)
\psline{->}(2.6,5)(3.3,5)
\psline{->}(4.6,5)(5.3,5)
\psline{->}(6.6,5)(7.3,5)

\psline(1.95,4.6)(1.95,4)
\psline(2.1,4.6)(2.1,4)
\psline{->}(4,4.6)(4,4)
\psline{->}(6,4.6)(6,4)

\uput[l](0.5,3.5){0}
\uput[l](2.5,3.5){$X$}
\uput[l](4.5,3.5){$M$}
\uput[l](6.5,3.5){$Y$}
\uput[l](8.15,3.5){$0$}
\psline{->}(0.6,3.5)(1.3,3.5)
\psline{->}(2.6,3.5)(3.3,3.5)
\psline{->}(4.6,3.5)(5.3,3.5)
\psline{->}(6.6,3.5)(7.3,3.5)

\psline{->}(4,3.1)(4,2.5)
\psline{->}(6,3.1)(6,2.5)

\uput[l](4.5,2){$Y''$}
\uput[l](6.5,2){$Y''$}
\psline(4.6,2.1)(5.3,2.1)
\psline(4.6,1.95)(5.3,1.95)

\psline{->}(4,1.6)(4,1)
\psline{->}(6,1.6)(6,1)
\uput[d](4,1){0}
\uput[d](6,1){0}
\end{pspicture}
\end{center}}\noindent
A similar result holds for an exact sequence
$0\to X'\to X\to X''\to 0$.

\end{lem}

\section{Nilpotent representations of a cyclic quiver}

Let $\dt=\dt_n$ ($n\geq 2$) be the cyclic quiver with $n$ vertices:

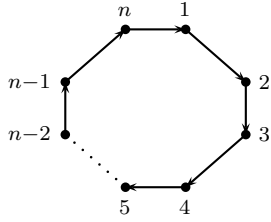
\begin{center}
\begin{pspicture}(0,-.6)(2,2.6)
\psset{xunit=.8cm,yunit=.7cm}
\psdot*(1,0)
\psdot*(2,0)
\psdot*(3,1)
\psdot*(3,2)
\psdot*(2,3)
\psdot*(1,3)
\psdot*(0,2)
\psdot*(0,1)
\uput[l](0,1){$_{n-2}$}
\uput[l](0,2){$_{n-1}$}
\uput[u](1,3){$_n$}
\uput[u](2,3){$_1$}
\uput[r](3,2){$_2$}
\uput[r](3,1){$_3$}
\uput[d](2,0){$_4$}
\uput[d](1,0){$_5$}
\psline{->}(1,3)(2,3)
\psline{->}(2,3)(3,2)
\psline{->}(3,2)(3,1)
\psline{->}(3,1)(2,0)
\psline{->}(2,0)(1,0)
\psline[linestyle=dotted,linewidth=1pt](1,0)(0,1)
\psline{->}(0,1)(0,2)
\psline{->}(0,2)(1,3)
\end{pspicture}
\end{center}

\noindent
Then $\dt$ determines the Cartan datum $(I,\centerdot)$ with
$I=\{1,2,\cdots,n\}$ and $i\centerdot i=2$ and $2\frac{i\centerdot j}
{i\centerdot i}=-\delta_{j,i+1}$ ($i,j\in I$).
It is always understood that $I$ is identified with
$\bbz/n\bbz=\{\bar 1,\bar 2,\cdots,\bar n\}$. When no risk arises,
we always drop the bars on the elements of $I$ for notational simplicity.
Thus, $n+1$ is understood as 1.

Let $k$ be a field. A representation $V=(V_i,f_i)_i$ of $\dt$ over $k$
is called {\it nilpotent} if the composition $f_n\cdots f_2f_1:V_1\ra V_1$
is nilpotent or equivalently, one of the
$f_{i-1}\cdots f_nf_1\cdots f_i:V_i\ra V_i$, $2\leq i\leq n$,
is nilpotent. The dimension vector of $V$ is defined to be $\udim V=
(\dim_k V_i)_i\in\bbn^n$. By $\bbt=\bbt(n,k)$ we denote the category of
finite-dimensional nilpotent representations of $\dt$ over $k$. Then $\bbt$ is an
abelian subcategory of $\rep\dt$, the category of all finite-dimensional
representations of $\dt$. Thus, the objects in $\bbt$ are also called modules.
Note that each nilpotent representation $M$ with dimension
$d$ can be considered as a module over finite dimensional self-injective algebras
$k\dt/J^m$ for $m\geq d+1$, where $J$ denotes the ideal of the path algebra $k\dt$
generated by all arrows of $\dt$.

For each vertex $i\in I$, we have a one-dimensional representation (or simple
module) $S_i=S_{ik}$. These $S_i$ form a complete set of simple objects in $\bbt$.
Further, for each integer $l\geq 1$, there is a unique
(up to isomorphism) indecomposable representation $S_i[l]$ in $\bbt$ of length $l$
with top $S_i$. It is well known that $S_i[l]$, $i\in I, l\geq 1$, yield all
isoclasses of indecomposable representations in $\bbt$. Then all isoclasses of
representations in $\bbt$ can be indexed by the set of $n$-tuples of partitions
which we define in the following.

A partition is a sequence $p=(p_1,p_2,\cdots,p_m,\cdots)$ of non-negative integers
in decreasing order
$$p_1\geq p_2\geq\cdots\geq p_m\geq\cdots$$
and containing only finitely many non-zero terms. The non-zero $p_i$
are called the {\it parts} of $p$, and we say that the $i$th part is of {\it length} $p_i$.
If $p_{m+1}=0$, we will
usually write $p=(p_1,p_2,\cdots,p_m)$. Thus we identify finite sequences which
only differ by adding some zeros at the end. Note that  the empty partition
$\emptyset:=(0,0,\cdots)$ is the only partition with $0$ part. By $\tilde p$ we denote the dual partition of $p$.
Finally, we denote by $\Pi$ the set of $n$-tuples
of partitions. Then each $n$-tuple $\pi=(\pi^{(1)}, \pi^{(2)},\cdots, \pi^{(n)})$ of
partitions defines a representation
$$M(\pi)=M_k(\pi)=\bigoplus_{{i\in I}\atop{j\geq 1}}  S_i[{\tilde \pi}^{(i)}_j]$$
in $\bbt(n,k)$, where $\tpi^{(i)}=(\tpi^{(i)}_1,\tpi^{(i)}_2,\cdots)$ is the dual
partition\footnote{Note that we adopt the notation from \cite[3.3]{R93} in order to view a module as
a column in the Young diagram of a partition.} of $\pi^{(i)}$, $i\in I$.
In this way we obtain a bijection between $\Pi$ and the set
of isoclasses of representations in $\bbt$. Note that this bijection is
independent of the field $k$.

For each module $M$ in $\bbt$, we denote by $[M]$ the isoclass of
$M$. For $a\geq 1$, we write
$$aM:=\underbrace{M\oplus\cdots\oplus M}_a.$$

For each partition $p=(p_1,p_2,\cdots,p_m)$ and each
$i\in I$, we set
$$M_i(p)=\bigoplus_{j=1}^m S_i[p_j].$$
Thus, if $\pi=(\pi^{(1)}, \pi^{(2)},\cdots, \pi^{(n)})\in\Pi$, we have
$$M(\pi)=\bigoplus_{i=1}^n M_i(\tpi^{(i)}).$$
The following easy facts about nilpotent representations turn out
to be very useful. For the convenience of the reader,
we include some proof.

\begin{lem}\label{EXT} Let $p=(p_1,p_2,\cdots,p_m)$ be a partition with $m$
parts and assume $i\in I$.

(1) A representation $L$ is an extension of $S_i$ by $M_{i+1}(p)$ if and
only if
$$L\cong L_t:=S_i[p_t+1]\oplus\bigoplus_{j:j\not=t} S_{i+1}[p_j]$$
for some $1\leq t\leq m+1$, where $p_{m+1}=0$.

(2) A submodule $N$ of $M_i(p)$ is maximal if and only if
$$N\cong N_t:= S_{i+1}[p_t-1]\oplus\bigoplus_{j:j\not=t} S_i[p_j]$$
for some $1\leq t\leq m$.

(3) If we choose $v_j\in S_i[p_j]\backslash\rad(S_i[p_j])$ and
$x_1,\cdots, x_m\in k$ are not all zero, then the element
$x_1v_1+\cdots+ x_mv_m\in M_i(p)$ generates a submodule isomorphic to
$S_i[p_a]$, where $a$ is the smallest index such that $x_a\not=0$.
\end{lem}

\begin{pf} (1)  Clearly, each $L_t$ admits an exact sequence
$$0\lra M_{i+1}(p)\lra L_t\lra S_i\lra 0,$$
that is, $L_t$ is an extension of $S_i$ by $M_{i+1}(p)$.
Conversely, let $L$ be any extension of $S_i$ by $M_{i+1}(p)$. We
may consider $M_{i+1}(p)$ as a submodule of $L$. Then
$$\rad (M_{i+1}(p))\subseteq \rad (L)\subseteq M_{i+1}(p).$$
It follows that
$$\top (L)=L/\rad (L)\cong S_i\oplus M_{i+1}(p)/\rad (L)
\cong S_i\oplus a S_{i+1}$$ for some $a\geq 0$. Applying
Krull-Remak-Schmidt theorem to $L$, we conclude that $L$ is
isomorphic to some $L_t$.

(2) If $N\cong N_t$ for some $1\leq t\leq m$, it is obvious that
$N$ is a maximal submodule of $M:=M_i(p)$. Conversely, let $N$ be a maximal
submodule of $M$. Then $\rad (M)\subseteq N$ and there are exact sequences
$$0\lra  N/\rad (M)\lra M/\rad (M)\lra M/N\cong S_i\lra 0$$
and
$$0\lra  \rad (M)/\rad (N)\lra N/\rad (N)\lra N/\rad (M)\lra 0.$$
It follows that $N/\rad (N)\cong (m-1)S_i\oplus aS_{i+1}$ for some $a\geq 0$.
The statement (1) implies $a\leq 1$ since $M$ is an extension of $S_i$
by $N$. Hence, $N$ is isomorphic to some $N_t$.

(3) This is obvious.
\end{pf}

The next lemma describes extensions of indecomposable nilpotent
representations by indecomposable ones.

\begin{lem}\label{INEXT} Let $i,j\in I$ and $l,m\geq 1$. If
there is a non-split exact sequence
\begin{equation}\label{split}
0\lra S_j[m]\lra M\lra S_i[l]\lra 0,
\end{equation}
then there are integers $r\geq 1, s\geq 0$, and $t\geq 1$
such that $j={\overline {i+r}}, l=r+s, m=s+t$ and
$M\cong S_i[r+s+t]\oplus S_{\overline {i+r}}[s].$
\end{lem}

\begin{pf} We apply induction on $l$. If $l=1$, then $j=\overline {i+1}$.
We take $r=1, s=0$, and $t=m$, as required. Let now $l>1$. Consider
the following exact sequence
$$0\lra S_{\overline {i+1}}[l-1]\lra S_i[l]\lra S_i\lra 0.$$
This together with (\ref{split}) induces by Lemma \ref{pbk}
exact sequences
\begin{align}\label{splita}
0\ra S_j[m]\ra &M'\ra S_{\overline {i+1}}[l-1]\ra 0\\
0\ra &M'\ra M\ra S_i\ra0.\label{splitb}
\end{align}
If (\ref{splita}) splits
and $j\not=\overline {i+1}$, then (\ref{split}) splits, contrary to
the assumption. Therefore, $j=\overline {i+1}$, $m\geq l$,
and $M\cong S_i[m+1]\oplus S_{\overline {i+1}}[l-1]$.
Thus, we can take $r=1, s=l-1$ and $t=m-l+1$.
If (\ref{splita}) does not split,
then, by induction, there are integers $r'\geq 1, s'\geq 0$, and
$t'\geq 1$ such that $j={\overline {i+1+r'}}, l-1=r'+s, m=s'+t'$ and
$$M'\cong S_{\overline {i+1}}[r'+s'+t']\oplus S_{\overline {i+1+r'}}[s'].$$
By Lemma \ref{EXT}(1), (\ref{splitb}) implies
$M\cong S_i[r'+s'+t'+1]\oplus S_{\overline {i+1+r'}}[s'].$
By putting $r=r'+1, s=s'$ and $t=t'$, we get
$$M\cong S_i[r+s+t]\oplus S_{\overline {i+r}}[s].$$
\end{pf}

Extensions (\ref{split}) can be used (cf. \cite[4.7]{R93}) to define
a relation $\dpeq$ on $\Pi$
by setting $\mu\dpeq\lz$
if there exist a $\pi\in\Pi$ and a non-split extension (\ref{split})
such that $M(\lz)\cong M(\pi)\oplus M$ and $M(\mu)\cong M(\pi)\oplus
S_j[m]\oplus S_i[l]$. Thus, by Lemma \ref{INEXT}, we see that
$\mu\dpeq\lz$ if and only if there are
integers $r\geq 1, s\geq 0, t\geq 1$, $i\in I$, and $\pi\in\Pi$
such that (a) $\pi$ is obtained from $\lz$ by deleting a part of length
$r+s+t$ in ${\tilde \lz}^{(i)}$ and a part of length $s$ in
${\tilde \lz}^{(\overline {i+r})}$, and (b) $\pi$ is also obtained from
$\mu$ by deleting a part of length $r+s$ in ${\tilde \mu}^{(i)}$ and
another part of length $s+t$ in ${\tilde \mu}^{(\overline {i+r})}$.
In other words, we have
$$\aligned
&M(\lz)=M(\pi)\oplus S_i[r+s+t]\oplus S_{\overline {i+r}}[s]\cr
\text{and}\;\;\;
&M(\mu)=M(\pi)\oplus S_i[r+s]\oplus S_{\overline {i+r}}[s+t].
\endaligned$$
In this case, there is an exact sequence
$$0\lra S_{\overline {i+r}}[s+t]\lra M(\lz)\lra M(\pi)\oplus S_i[r+s]\lra 0.$$
Let $\peq$ be the partial ordering on $\Pi$ generated by the relation
$\dpeq$. Thus, $\mu\prec\lz$ means that there are $\mu_0=\mu,\mu_1,
\cdots,\mu_m=\lz$ in $\Pi$ such that
$\mu_0\dpeq \mu_1\dpeq\cdots\dpeq\mu_m$.
We called $\peq$ the {\it extension order} on $\Pi$.

\section{Generic extensions and the monoid $\cm$}

In this section $k$ is algebraically closed. We study generic extensions of nilpotent representations.
Their existence follows from a similar argument as in \cite[Sect. 2]{Re}.

Let $\bd=(d_i)_i\in\bbn^n$. Then each representation $M=(k^{d_i},f_i)_i$
(not necessarily nilpotent) of $\dt$ can be identified with the point
 $$f=(f_i)_i\in\prod_{i=1}^n\hom_k(k^{d_i},k^{d_{i+1}})\cong
\prod_{i=1}^n k^{d_{i+1}\times d_i}=:R(\bd),$$
where $d_{n+1}=d_1$. The algebraic group
$$GL(\bd)=\prod_{i=1}^n GL_{d_i}(k)$$
acts on $R(\bd)$ by conjugation
$$(g_i)_if=(g_{i+1}f_ig_i^{-1})_i,$$
where $g_{n+1}=g_1$. Then the $GL(\bd)$-orbits in $R(\bd)$ correspond bijectively to
the isoclasses of representations of $\dt$ with dimension vector $\bd$.

The stabilizer $GL(\bd)_M=\{g\in GL(\bd)|gM=M\}$ of $M$ is the automorphism
group of $M$ which
is Zariski-open in $\mbox{End}(M)$ and has dimension equal to
$\dim\mbox{End}(M)$. It follows that
the orbit $\co_M$ of $M$ has dimension
$$\dim\co_M=\dim GL(\bd) -\dim\mbox{End}(M).$$
For two representations $M,N$, we say that $M$ degenerates to $N$,
or that $N$ is a {\it degeneration} of $M$, and write\footnote{The
order $\led$ here is opposite to the so-called degeneration order
used, e.g., in \cite{Re}. We follow a traditional notation for the
``Bruhat-Chevalley order'' of a Coxeter group.} $N\led M$, if
$\co_N\subseteq {\overline \co_M}$, the closure of $\co_M$. Note
that $N<M \iff \co_N\subseteq {\overline \co_M}\backslash\co_M$.
The following lemma is well-known (see for example
\cite[1.1]{Bo}).

\begin{lem} \label{BOD}If $0\to N\to E\to M\to 0$ is exact and non-split, then
$M\oplus N< E$.
\end{lem}

Since nilpotent representations with fixed dimension can be considered
as modules over a finite dimensional representation-finite self-injective
algebra $k\dt/J^m$, \cite[Thm 2]{Zwa} leads to the following
alternative descriptions  of degenerations  of nilpotent representations.

\begin{lem} \label{HOM} Let $M$ and $N$ be nilpotent representations.
Then the following statements are equivalent:

(1) $N\leq M$,

(2) $\dim_k\hom(N,X)\geq\dim_k\hom (M,X)$ for all $X\in\bbt$,

(3) $\dim_k\hom(X, N)\geq\dim_k\hom (X, M)$ for all $X\in\bbt$,

(4) there are representations $M_i,X_i,Y_i$ ($1\leq i\leq t$) in $\bbt$ and short exact
sequences $0\ra X_i\ra M_i\ra Y_i\ra 0$ such that $M=M_1,
M_{i+1}=X_i\oplus Y_i$, $1\leq i\leq t$, and $N=M_{t+1}$.
\end{lem}

We remark that,
since the dimension $\dim_k\hom(X, Y)$
is independent of $k$  for any $X,Y\in\bbt(n,k)$, the lemma guarantees that the ordering $\led$
 on nilpotent
representations is well-defined over any field $k$.

Further, let $R_0(\bd)$ be the subset of $R(\bd)$ consisting of nilpotent representations.
It is easy to see that $R_0(\bd)$ is a closed subvariety of $R(\bd)$
stable under
$GL(\bd)$ and contains only finitely many orbits.
Clearly, an extension of a nilpotent representation by a nilpotent
one is nilpotent. Hence,
we have immediately the following (cf. \cite[2.1,2.3]{Re}). We refer their
proofs to \cite{Re}.

\begin{lem}\label{GEL} For $\bd,\bd',\bd''\in{\Bbb N}^n$ with $\bd=\bd'+\bd''$ and
subsets $\cx\subset R_0(\bd')$ and $\cy\subset R_0(\bd'')$ stable under
$GL(\bd')$ and $GL(\bd'')$, respectively, let $\ce=\ce(\cx,\cy)$ be the subset
of $R_0(\bd)$ of all extensions of representations $M\in\cx$ by
representations $N\in\cy$. Then $\ce$ is $GL(\bd)$ stable, and hence,
constructible.  Moreover, if $\cx$ and $\cy$ are irreducible (resp. closed),
then so is $\ce$.
\end{lem}

Applying the lemma to $\cx=\co_M$ and $\cy=\co_N$, we see that,
for any two nilpotent representations $M$ and $N$, there exists
a unique extension $G$ (up to isomorphism)  of $M$ by $N$
with maximal dimension of  $\co_G$, or equivalently,
with minimal $\dim\text{End}(G)$.
This representation $G$ is called a {\it generic extension} of
$M$ by $N$ and is denoted
by $M\ast N:=G$.

It can be seen easily that
 \cite[2.4]{Re} continues to hold in this case.

\begin{prop} \label{GE}
Let $M,N,X$ be nilpotent representations of $\dt$. Then
$$X\led M\ast N \iff X\in\ce(\co_{M'},\co_{N'})\text{ for some }  M'\led M,N'\led N.$$
In particular, we have
$$M'\led M,N'\led N\Longrightarrow M'*N'\led M*N.$$
\end{prop}

As a conclusion of Lemma \ref{GEL} and Prop. \ref{GE}, we can define a monoid
$\cm=\cm(\dt)$ whose elements
are the isoclasses of nilpotent representations with multiplication $[M]\ast[N]=[M\ast N]$
and unit element $1_\cm=[0]$, the isoclass of the zero representation. Note that the
associativity of the multiplication $\ast$ follows from
Lemma \ref{pbk} and Prop. \ref{GE} (cf. \cite[3.1]{Re}). Inductively, we
see easily that for any nilpotent representations
$M_1,M_2,\cdots, M_t$, the representation $M_1\ast M_2\ast\cdots \ast M_t$ is the nilpotent
representation $G$ ( unique up to isomorphism) with minimal dimension of its endomorphism
algebra such that $G$ has a filtration
$$G=G_0\supset G_1\supset \cdots\supset G_{t-1}\supset G_t=0$$
with $G_{s-1}/G_s\cong M_s$, $1\leq s\leq t$.

Unlike the situation in \cite{Re}, in our case the submonoid
$\cm_c$ of $\cm$ generated by $[S_i]$, $i\in I$, is proper. However, we have the
following similar result to \cite[3.4]{Re}.

\begin{prop}\label{reln} (1) In case $n\geq 3$, the following relations hold in $\cm_c$:
$$\aligned
&[S_i]\ast [S_j]=[S_j]\ast [S_i]\;\;\mbox{if $j\not\equiv i\pm 1$ (mod $n$),
$i,j\in I$},\cr
&[S_i]\ast [S_{i+1}]\ast [S_i]=[S_i]\ast [S_i]\ast [S_{i+1}],\cr
&[S_{i+1}]\ast [S_i]\ast [S_{i+1}]=[S_i]\ast [S_{i+1}]\ast [S_{i+1}],\;
i\in I.
\endaligned$$

(2) In case $n=2$, the following relations hold in $\cm_c$:
$$\aligned
&[S_1]\ast [S_2]\ast [S_1]\ast [S_1]=[S_1]\ast [S_1]\ast [S_2]\ast [S_1],\cr
&[S_2]\ast [S_1]\ast [S_2]\ast [S_2]=[S_2]\ast [S_2]\ast [S_1]\ast [S_2].
\endaligned$$
\end{prop}

The proof of this proposition is straightforward.

Our next purpose is to describe the generic extension of a simple representation
by any nilpotent representation.

For each $i\in I$, we define a map
$$\sz^+_i:\Pi\ra\Pi$$
as follows: for
$\pi=(\pi^{(1)},\pi^{(2)},\cdots,\pi^{(n)})\in\Pi$, $\sz^+_i\pi=:\lambda$ is
defined by
$$
\lambda^{(j)}=\cases \pi^{(j)}&\text{ if }j\neq i,i+1\cr
\pi^{(i)}+{\mathbf 1}_{\tpi^{(i+1)}_1+1}&\text{ if }j=i\cr
\pi^{(i+1)}-{\mathbf 1}_{\tpi^{(i+1)}_1}&\text{ if }j=i+1,\cr
\endcases
$$
where ${\mathbf 1}_a=(\underbrace{1,\cdots,1}_a).$
Note that, if we define $\lz$ dually, then
$\tilde\lambda^{(i+1)}$ is
obtained
by deleting the part $\tpi^{(i+1)}_1$ from $\tpi^{(i+1)}$, and
$\tilde\lambda^{(i)}$ is
obtained by adding the part
$\tpi^{(i+1)}_1+1$ to $\tpi^{(i)}$.

\begin{example}\label{ALG} Let $n=3$ and $\pi=(\pi^{(1)},\pi^{(2)}, \pi^{(3)})$ with
$$\pi^{(1)}=(4,3,3,1,1),\, \pi^{(2)}=(3,2,1),\,\pi^{(3)}=(2,2).$$
It is easy to see that  $\tpi^{(1)}_1=5, \tpi^{(2)}_1=3, \tpi^{(3)}_1=2$.
Then
$$\aligned
\sz^+_1\pi&=\bigl((5,4,4,2,1),\, (2,1),\,(2,2)\bigr),\cr
\sz^+_2\pi&=\bigl((4,3,3,1,1),\, (4,3, 2),\,(1,1)\bigr),\cr
\sz^+_3\pi&=\bigl((3,2,2),\, (3, 2,1),\,(3,3,1,1,1,1)\bigr).\cr
\endaligned$$
Intuitively, for example, $\pi$ and $\sz_1^+\pi$ are illustrated by
their Young diagrams as follows:

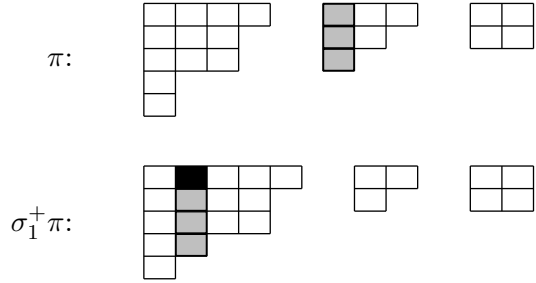
\begin{center}
\begin{pspicture}(0,0)(6,2)
\psset{xunit=.7cm,yunit=.5cm,linewidth=0.5pt}
\uput[l](0.1,1.5){$\pi$:}

\psline(1.2,3)(3.6,3)
\psline(1.2,2.4)(3.6,2.4)
\psline(1.2,1.8)(3,1.8)
\psline(1.2,1.2)(3,1.2)
\psline(1.2,0.6)(1.8,0.6)
\psline(1.2,0)(1.8,0)
\psline(1.2,3)(1.2,0)
\psline(1.8,3)(1.8,0)
\psline(2.4,3)(2.4,1.2)
\psline(3,3)(3,1.2)
\psline(3.6,3)(3.6,2.4)

\psline(4.6,3)(6.4,3)
\psline(4.6,2.4)(6.4,2.4)
\psline(4.6,1.8)(5.8,1.8)
\psline(4.6,1.2)(5.2,1.2)
\psline(4.6,3)(4.6,1.2)
\psline(5.2,3)(5.2,1.2)
\psline(5.8,3)(5.8,1.8)
\psline(6.4,3)(6.4,2.4)
\psframe[fillstyle=solid,fillcolor=lightgray](4.6,1.2)(5.2,1.8)
\psframe[fillstyle=solid,fillcolor=lightgray](4.6,1.8)(5.2,2.4)
\psframe[fillstyle=solid,fillcolor=lightgray](4.6,2.4)(5.2,3)

\psline(7.4,3)(8.6,3)
\psline(7.4,2.4)(8.6,2.4)
\psline(7.4,1.8)(8.6,1.8)
\psline(7.4,3)(7.4,1.8)
\psline(8,3)(8,1.8)
\psline(8.6,3)(8.6,1.8)

\end{pspicture}
\end{center}

\begin{center}
\begin{pspicture}(0,-0.2)(6,2)
\psset{xunit=.7cm,yunit=.5cm,linewidth=0.5pt}
\uput[l](0.1,1.5){$\sz_1^+\pi$:}

\psline(1.2,3)(4.2,3)
\psline(1.2,2.4)(4.2,2.4)
\psline(1.2,1.8)(3.6,1.8)
\psline(1.2,1.2)(3.6,1.2)
\psline(1.2,0.6)(2.4,0.6)
\psline(1.2,0)(1.8,0)
\psline(1.2,3)(1.2,0)
\psline(1.8,3)(1.8,0)
\psline(2.4,3)(2.4,0.6)
\psline(3,3)(3,1.2)
\psline(3.6,3)(3.6,1.2)
\psline(4.2,3)(4.2,2.4)
\psframe[fillstyle=solid,fillcolor=lightgray](1.8,0.6)(2.4,1.2)
\psframe[fillstyle=solid,fillcolor=lightgray](1.8,1.2)(2.4,1.8)
\psframe[fillstyle=solid,fillcolor=lightgray](1.8,1.8)(2.4,2.4)
\psframe*(1.8,2.4)(2.4,3)

\psline(5.2,3)(6.4,3)
\psline(5.2,2.4)(6.4,2.4)
\psline(5.2,1.8)(5.8,1.8)
\psline(5.2,3)(5.2,1.8)
\psline(5.8,3)(5.8,1.8)
\psline(6.4,3)(6.4,2.4)

\psline(7.4,3)(8.6,3)
\psline(7.4,2.4)(8.6,2.4)
\psline(7.4,1.8)(8.6,1.8)
\psline(7.4,3)(7.4,1.8)
\psline(8,3)(8,1.8)
\psline(8.6,3)(8.6,1.8)

\end{pspicture}
\end{center}

\noindent Here the coloured column in $\sz_1^+\pi$ is obtained from the
coloured one in $\pi$ topped with the black box.
\end{example}

\begin{prop}\label{SGE}  Let $i\in I$ and $\pi\in\Pi$. Then
$$S_i\ast M(\pi)\cong M(\sz^+_i\pi).$$
\end{prop}

\begin{pf} Let $\pi=(\pi^{(1)}, \pi^{(2)},\cdots, \pi^{(n)})$. For each $a\in I$,
we set
$$M_a(\pi):=M_a(\tpi^{(a)})=\bigoplus_{j\geq 1} S_a[\tpi^{(a)}_j],$$
where $\tpi^{(a)}=(\tpi^{(a)}_1,\tpi^{(a)}_2,\cdots)$ is the dual partition of
$\pi^{(a)}$. Then we have $M(\pi)=\oplus_{a=1}^n M_a(\pi)$. Since
$\ext^1(M_a(\pi),S_i)=0$ for all $a\not=i+1$, it follows that
$$S_i\ast M(\pi)\cong
\bigoplus_{a\not=i+1}M_a(\pi)\oplus (S_i\ast M_{i+1}(\pi)).$$
To prove the proposition, it suffices to prove
\begin{equation}\label{SGEa}
S_i\ast M_{i+1}(\pi)\cong S_i[\tpi^{(i+1)}_1+1]\oplus
\bigoplus_{j\geq 2} S_{i+1}[\tpi^{(i+1)}_j].
\end{equation}
If $\pi^{(i+1)}=\emptyset$, this is obvious. Let now $\pi^{(i+1)}\not=\emptyset$.
For simplicity, we write $\tpi^{(i+1)}=(p_1,p_2\cdots,p_m)$ with
$p_1\geq p_2\cdots\geq p_m>p_{m+1}=0$, where $p_t=\tpi^{(i+1)}_t$, $1\leq t\leq m$.
Then, by Lemma \ref{EXT}(1), each extension of $S_i$ by $M_{i+1}(\pi)$ is isomorphic to
$$L_t=S_i[p_t+1]\oplus\bigoplus_{j\not=t} S_{i+1}[p_j]$$
for some $t$ with $1\leq t\leq m+1$.

Suppose $1\leq t\leq m-1$. If $p_t=p_{t+1}$, it is obvious
that $L_t\cong L_{t+1}$.
If $p_t>p_{t+1}$, we get a non-split exact sequence
$$0\lra S_{i+1}[p_t]\lra S_i[p_t+1]\oplus S_{i+1}[p_{t+1}]
\lra S_i[p_{t+1}+1]\lra 0.$$
This yields a non-split exact sequence
$$0\lra \bigoplus_{j\not=t+1} S_{i+1}[p_j]\lra L_t\lra S_i[p_{t+1}+1]\lra 0.$$
Hence, by Lemma \ref{BOD}
$$L_{t+1}=S_i[p_{t+1}+1]\oplus\bigoplus_{j\not=t+1} S_{i+1}[p_j]<L_t.$$
Consequently, we have
$$S_i\oplus M_{i+1}(\pi)=L_{m+1}\led L_m\led \cdots\led L_2\led L_1\cong
S_i\ast M_{i+1}(\pi),$$
proving (\ref{SGEa}).
\end{pf}

\begin{rem} \label{GEk} (1) Dually, one can describe the generic extension
of any nilpotent representation
by a simple representation. We leave the detail to the reader.

(2) We observe from Lemma \ref{EXT} that the
number of non-isomorphic extensions of a simple module by any given
nilpotent representation is independent
of the field $k$.  This together with the fact that
$\dim\text{End}(M\oplus N)>\text{End}(E)$ for a non-split exact
sequence $0\to N\to E\to M\to 0$ in $\bbt$ implies from
the proof of \ref{SGE} that, for any field $k$ and $M_k\in\bbt(n,k)$,
$$(S_{ik}*M_k)\otimes \bar k\cong S_{i\bar k}*M_{\bar k},$$
where $S_{ik}*M_k$ is the extension of $S_{ik}$ by $M_k$ with minimal
dimension of $\text{End}(S_{ik}*M_k)$.
\end{rem}

\section{The submonoid $\cm_c$}

Following \cite[4.1]{R93}, an $n$-tuple $\pi=(\pi^{(1)},\pi^{(2)},\cdots,\pi^{(n)})$
of partitions
is called {\it separated}
if, for each $t\geq 1$, there is some $a(t)=a(t,\pi)\in I$ such that
$\tpi^{(a(t))}_j\neq t$ for all $j\geq 1$. By $\Pi^s$ we denote the set of separated
$n$-tuples of partitions. Then each
$\sz^+_i$, $i\in I$ defined in \S 2 induces a map
$$\sz^+_i:\Pi^s\ra\Pi^s.$$
Indeed,
let $\pi\in\Pi^s$, for each $t\geq 1$, if $t\not=\tpi^{(i+1)}_1+1$,
$({\widetilde {\sz^+_i\pi}})^{(a(t,\pi))}$  has no part equal to $t$, and if
$t=\tpi^{(i+1)}_1+1$, $({\widetilde {\sz^+_i\pi}})^{(i+1)}$ has no part
equal to $t$, proving $\sz^+_i\pi\in\Pi^s$.

We observe that, for $\pi=(\pi^{(1)},\pi^{(2)},\cdots,\pi^{(n)})\in\Pi^s$,
if $\pi\not=\emptyset$, the $n$-tuple of empty partitions,
 then $\tpi^{(i)}_1$, $i\in I$, can not be the same,
and there exists
an $i\in I$ (not necessarily unique) such that $\tpi^{(i)}_1>\tpi^{(i+1)}_1$.
(Note that $\tpi^{(n+1)}=\tpi^{(1)}$ and that $\tpi^{(i+1)}_1=0$ if
$\tpi^{(i+1)}$ is the empty partition.) Thus, if we put
$\Pi^s_i=\sz^+_i(\Pi^s)$, we have immediately from the definition that
$$\Pi^s=\{\emptyset\}\cup\bigcup_{i=1}^n\Pi^s_i.$$

As in \cite[2.1]{R93}, let $\ooz$ denote the free semigroup with unit element,
generated by the set $I=\{1,2,\cdots,n\}$. The elements in $\ooz$ are called
{\it words}
and are of the form $w=i_1i_2\cdots i_m$ with $i_1,i_2,\cdots,i_m\in I$
and $m\geq 0$. For each $w=i_1i_2\cdots i_m\in\ooz$, we set
$$M(w)=S_{i_1}\ast S_{i_2}\ast\cdots\ast S_{i_m}.$$
Then there is a unique $\pi\in\Pi$ such that $M(w)\cong M(\pi)$, and we set
$\wp(w)=\pi$. In this way we obtain a map
$$\wp:\ooz\to \Pi, w\mapsto\pi=\wp(w).$$
Note that by Remark \ref{GEk}(2) the map $\wp$ is independent of the field
$k$.

\begin{thm}\label{wp} (1) The map $\wp$ induces a surjection
$$\wp:\ooz\twoheadrightarrow\Pi^s.$$

(2) Let $\cm(\Pi^s):=\{[M(\pi)]:\pi\in\Pi^s\}.$ Then $\cm_c=\cm(\Pi^s)$.
\end{thm}

\begin{pf} By Prop. \ref{SGE} and induction on $m$, we see easily that,
for $w=i_1i_2\cdots i_m\in\ooz$,
\begin{equation}
\wp(w)=\sz^+_{i_1}\sz^+_{i_2}\cdots\sz^+_{i_m}(\emptyset)\in\Pi^s,
\end{equation}
since $\emptyset$
is separated. Consequently,
we obtain the inclusion $\cm_c\subseteq\cm(\Pi^s)$.
Clearly, if $\wp$ maps onto $\Pi^s$, then the inclusion
$\cm(\Pi^s)\subseteq \cm_c$ follows.
It remains to prove that $\wp$ is surjective.

The surjective map $\sz_i^+:\Pi^s\to\Pi^s_i$ has a right inverse
$\sz_i^-:\Pi^s_i\to \Pi^s$ which can be defined as follow: for
$\pi\in \Pi^s_i$,  define $\mu:=\sz_i^-(\pi)\in\Pi$ by
\begin{equation}
\mu^{(a)}=\cases \mu^{(a)}&\text{ if }a\neq i,i+1\cr
\mu^{(i)}-{\mathbf 1}_{\tpi^{(i)}_j}&\text{ if }a= i\cr
\mu^{(i+1)}+{\mathbf 1}_{\tpi^{(i)}_j-1}&\text{ if }a= i+1\cr
\endcases
\end{equation}
where $\tpi^{(i)}_j$ is the smallest part in $\tpi^{(i)}$
such that $\tpi^{(i)}_j>\tpi^{(i+1)}_1$.
(If we define $\mu$ dually, then
$\tilde \mu^{(i)}$ is the partition obtained from
$\tpi^{(i)}$ by deleting the part
$\tpi^{(i)}_j$ and $\tilde \mu^{(i+1)}$ is obtained by adding the part
 $\tpi^{(i)}_j-1$ to $\tpi^{(i+1)}$.)  In other words,
the corresponding module $M(\mu)$ is obtained
by deleting the simple top $S_i$ of the
indecomposable summand of $M(\pi)$ corresponding to the $j$th column
$\tpi^{(i)}_j$ of $\tpi^{(i)}$.
It is easy to see that $\mu$ is separated and $\sz^+_i(\sz_i^-\pi)=\pi$.

We now can associate to each $\pi\in\Pi^s$ a word $w$
with $\wp(w)=\pi$ (and so $m=|\pi|:=\sum_{i\in I,j\geq 1}\pi^{(i)}_j$).
If $\pi=\emptyset$, we set $w=1$. Otherwise, $\pi\in\Pi^s_{i_1}$
for some $i_1$ with
$\sz_{i_1}^-(\pi)\in\Pi_{i_2}^s$, etc. So we obtain
a sequence of $i_1,\cdots,
i_m$ such that
$$\sz_{i_m}^-\cdots\sz_{i_1}^-(\pi)=\emptyset.$$
Let $w=i_1\cdots i_m\in\ooz$. Then $\wp(w)=\pi$. Hence, $\wp$ maps
onto $\Pi^s$.
\end{pf}

Note that, in general, $\sz_i^-(\sz^+_i\pi)\not=\pi$. For example,
let $n=2$ and
$$\pi=\bigl((\pi^{(1)}=(2,1,1),\,\pi^{(2)}=(1,1) \bigr).$$
Then $\sz_1^-(\sz^+_1\pi)=\bigl((2,2,2),\,\emptyset \bigr)$
 does not equal to $\pi$.

\begin{example}\label{WORD} Let us consider Example \ref{ALG}, that is, $n=3$ and
$$\pi=\bigl(\pi^{(1)}=(4,3,3,1,1),\,\pi^{(2)}=(3,2,1),\, \pi^{(3)}=(2,2)\bigr).$$
Since $\tpi^{(1)}_1=5, \tpi^{(2)}_1=3, \tpi^{(3)}_1=2$,
we have both $\pi\in\Pi^s_1$ and $\pi\in\Pi^s_2$. Then
$$\aligned
&\sz^-_1\pi=\bigl((3,2,2),\,(4,3,2,1)\,(2,2)\bigr),\cr
&\sz^-_2\pi=\bigl((4,3,3,1,1),\,(2,1)\,(3,3)\bigr).
\endaligned$$
Repeating the above process, we totally get nine words in this way
$$\begin{array}{rlll}
&12^21^23^42^41^53^32, &21^323^42^41^53^32, &21^3323^32^41^53^32,\\
&21^33^223^22^41^53^32, &21^33^32^531^53^32, &21^33^32^5131^43^32,\\
&21^33^32^51^231^33^32, &21^33^32^51^331^23^32, &21^33^32^51^43^412.
\end{array}$$
 However, they are not
all the words in the fibre
$\wp^{-1}(\pi)$. In fact, $\wp^{-1}(\pi)$ consists of 141 words.
\end{example}

\begin{rem} In \cite[4.4]{R93} Ringel considers a proper subset $\ooz^c$ of
$\ooz$ consisting of  condensed words (see the definition in
\cite[4.4]{R93}) and defines a map $\vez:\ooz^c\ra\Pi^s$. His map
$\vez$ is in fact the restriction of our map $\wp:\ooz\ra\Pi^s$.
This can be seen from \cite[3.4,5.1,Thm C]{R93} and results
\ref{order} and \ref{DOMI} below.
\end{rem}

\section{Distinguished words}

For each word $w=i_1i_2\cdots i_m\in\ooz$ and each
representation $M$ in $\bbt(n,k)$, a composition series
$$M=M_0\supset M_1\supset \cdots\supset M_{m-1}\supset M_m=0$$
of $M$ is said to be of type $w$ if $M_{j-1}/M_j\cong S_{i_j}$ for all
$1\leq j\leq m$. We denote by $\lr{w,M}$
the number of composition series of $M$ of type $w$. Note that
if $k$ is an infinite field, $\lr{w,M}$ may be infinite.

\begin{lem}\label{COM} Let $F$ be another filed. Then for each $w\in\ooz$
and each $\lz\in\Pi$, we have
$$\lr{w,M_k(\lz)}\not=0\Longleftrightarrow\lr{w,M_F(\lz)}\not=0.$$
\end{lem}

\begin{pf} Let $w=i_1i_2\cdots i_m$ and $\lr{w,M_k(\lz)}\not=0$.
We use induction on $m$ to show $\lr{w,M_F(\lz)}\not=0$.

If $m=0$ or $1$, this is obvious. Let $m>1$. Since $\lr{w,M_k(\lz)}\not=0$,
$M_k(\lz)$ admits a submodule $N$ such that $M_k(\lz)/N\cong S_{i_1k}$
and $\lr{w', N}\not =0$, where $w'=i_2\cdots i_m$. Then there is a unique
$\mu\in\Pi$ with $N\cong M_k(\mu)$. Thus, we obtain an exact sequence
$$0\lra M_k(\mu)\lra M_k(\lz)\lra S_{i_1k}\lra 0.$$
Now Lemma \ref{EXT}(2) implies that there exists a similar exact
sequence in $\bbt(n,F)$
$$0\lra M_F(\mu)\lra M_F(\lz)\lra S_{i_1F}\lra 0.$$
By induction,  $\lr{w',M_k(\mu)}=\lr{w',N}\not=0$
implies $\lr{w',M_F(\mu)}\not=0$. Consequently, $\lr{w,M_F(\lz)}\not=0$.

The converse is obtained by reversing the roles of $k$ and $F$.
\end{pf}

We recall from \cite[2.3]{R93} the definition of a reduced
filtration. Let $w=i_1i_2\cdots i_m$ be a word in $\ooz$. Then $w$
can be uniquely expressed in the {\it tight form}
$w=j_1^{e_1}j_2^{e_2}\cdots j_t^{e_t}$, where $e_r\geq 1$
$\forall\, r$, and $j_r\not=j_{r+1}$ for $1\leq r\leq t-1$. A
filtration
$$M=M_0\supset M_1\supset \cdots \supset M_{t-1}\supset M_t=0$$
of a nilpotent representation $M$ is called a {\it reduced} filtration
of type $w$ if
$M_{r-1}/M_r\cong e_rS_{j_r}$
for all $1\leq r\leq t$. Note that any reduced filtration of $M$ of type $w$
can be refined to a composition series of $M$ of type $w$. Conversely, given a
composition series of $M$ of type $w$, there is a unique reduced filtration
of $M$ of type $w$ such that the given composition series is a refinement
of this reduced filtration.

\begin{defn}\label{DISD} A word $w$ is called {\it distinguished}
 if the module $M_k(\wp(w))$ over an algebraically closed field $k$
has a unique reduced filtration of type $w$.
\end{defn}

\begin{rem} The condensed words constructed in \cite[4.4]{R93} are
distinguished words (see \cite[5.8]{R93}), but the converse is
not true. For example, let $n=3$ and
$$\pi=\bigl(\pi^{(1)}=(3,2,1),\,\pi^{(2)}=(1,1),\, \pi^{(3)}=(1)\bigr).$$
Then the fibre $\wp^{-1}(\pi)$ contains 7 words:
$$\begin{array}{rlllll}
&12^21^23^32, &21^33^22^23, &2121^23^32, &21^323^32, &21^33^2232,\\
&21^3323^22, &21^2213^32.  & & &
\end{array}$$
Among them, the first two words are the only condensed ones, but the first
five words are distinguished.
\end{rem}

In the following we are going to characterize all
distinguished words.
Recall from \S2 the notation $M_i(p)=\bigoplus_{j=1}^m S_i[p_j]$,
where $p$ is a partition and $i\in I$.
We start with the following lemma.

\begin{lem}\label{DISL} Let $p=(p_1,p_2,\cdots,p_m)$ be a partition with
$m$ parts. For a fixed subsequence $1\leq j_1< j_2<\cdots <j_t\leq m$, let
$$N=\bigoplus_{r=1}^t S_{i+1}[p_{j_r}-1]\oplus\bigoplus_{j\not=j_1,..,j_t} S_i[p_j].$$
Then $M_i(p)$ has a unique submodule isomorphic to $N$
if and only if $j_r=r$ for all $1\leq r\leq t$ and $p_t>p_{t+1}$.
Moreover, if this is the case and $k$ is algebraically closed, then
each submodule $N'$ of $M_i(p)$ with $M_i(p)/N'\cong tS_i$ satisfies $N\led N'$.
\end{lem}

\begin{pf}
Suppose that $j_r=r$ for all $1\leq r\leq t$ and $p_t>p_{t+1}$.
It is easy to see that, if $U_j$ denotes the unique
maximal submodule of $S_i[p_j]$, the submodule
$$X:=\bigoplus_{r=1}^t U_r\oplus\bigoplus_{j=t+1}^m S_i[p_j]$$
is isomorphic to $N$. Clearly, $M_i(p)/X\cong tS_i$. Let $X'$ be a submodule
of $M_i(p)$ isomorphic to $N$. Then there is a short exact sequence
$$0\lra X'\skrel{\iota}{\lra} M_i(p)\skrel{\phi}{\lra} tS_i\lra 0,$$
where $\iota$ denotes the inclusion.
Choose $v_j\in S_i[p_j]\backslash U_j$.
Since $\oplus_{r=1}^m U_r=\rad (M_i(p))\subseteq X'$,
Lemma \ref{EXT}(3) together with the condition $p_t>p_{t+1}$
implies  that $\phi(v_1),\cdots,
\phi(v_t)$ are linearly independent and that $\phi(v_j)=0$ for
$t+1\leq j\leq m$. This forces $X'=\Ker\phi=X$.

Conversely, suppose that $M$ has a unique submodule isomorphic to $N$
and suppose that
either $j_r>r$ for some $1\leq r\leq t$ or $j_r=r$ for all
$1\leq r\leq t$, but $p_t=p_{t+1}$.
In both cases, there exists an $l\not= j_r$ for $1\leq r\leq t$ such that
$p_l\geq p_{j_t}$. For each $x\in k$, we denote by $T_x$ the submodule
of $S_i[p_{j_t}]\oplus S_i[p_l]$ generated by  $U_{j_t}$ and $xv_t+v_l$.
It is easy to see that $T_x$ are pairwise distinct and all isomorphic to
$S_{i+1}[p_{j_t}-1]\oplus S_i[p_l]$ with
$(S_i[p_{j_t}]\oplus S_i[p_l])/T_x\cong S_i$. For each $x\in k$, we set
$$N_x=\bigoplus_{r=1}^{t-1} U_{j_r}\oplus T_x\oplus
\bigoplus_{{j\not=l,j_1,..,j_t}} S_i[p_j]$$
Then $N_x,x\in k,$ are pairwise distinct submodules of $M$ which are all
isomorphic to $N$. This is a contradiction. Thus, we have
$j_r=r$ for all $1\leq r\leq t$ and $p_t>p_{t+1}$.

To see the second assertion, let $k$ be an algebraically closed field,
and suppose
$j_r=r$ for all $1\leq r\leq t$ and $p_t>p_{t+1}$. Take any submodule $N'$
of $M_i(p)$ such that $M_i(p)/N'\cong tS_i$. By Lemma \ref{EXT}(2), there
are $1\leq l_1<l_2<\cdots<l_t\leq m$ such that
$$N'\cong \bigoplus_{r=1}^t S_{i+1}[p_{l_r}-1]\oplus\bigoplus_{j\not=l_1,..,l_t} S_i[p_j].$$
Let $s$ be the number of $l_r$'s with $l_r>t$. We use induction on $s$ to prove
$$\bigoplus_{r=1}^t S_{i+1}[p_r-1]\oplus\bigoplus_{j=t+1}^m S_i[p_j]=N\led N'.$$
If $s=0$, then $N'\cong N$. If $s\geq 1$, then $l_t\geq t+1$ and
there exists an $ a\leq t$
such that $a\not= l_r$ for all $1\leq r\leq t$ and $p_a\geq p_t>p_{l_t}$.
This inequality $p_a>p_{l_t}$
gives rise to a non-split exact sequence
$$0\lra S_{i+1}[p_a-1]\lra S_i[p_a]\oplus S_{i+1}[p_{l_t}-1]\lra S_i[p_{l_t}]\lra 0,$$
which implies by Lemma \ref{BOD}
$$S_i[p_{l_t}]\oplus S_{i+1}[p_a-1]<S_i[p_a]\oplus S_{i+1}[p_{l_t}-1].$$
Thus,
$$N'':=\bigoplus_{r=1}^{t-1} S_{i+1}[p_{l_r}-1]\oplus S_i[p_{l_t}]
\oplus S_{i+1}[p_a-1]\oplus\bigoplus_{j\not=a,l_1,..,l_t} S_i[p_j]<N'.$$
The inductive hypothesis finally shows $N\leq N''<N'$. This completes the proof.
\end{pf}

Let $w=i_1^{e_1}i_2^{e_2}\cdots i_t^{e_t}\in\ooz$ be
in the tight form.
For each $0\leq a\leq t$, we put
$w_a=i_{a+1}^{e_{a+1}}\cdots i_t^{e_t}\;\;\text{and}\;\;\nu(a)=\wp(w_a).$
In particular, $w_0=w$ and $w_t=1$. Further,  for $a\geq 1$, we have
$$\nu(a-1)=\wp(w_{a-1})=\underbrace{\sz^+_{i_a}\cdots \sz^+_{i_a}}_{e_a}
\nu(a).$$

\begin{thm}\label{DISP} Maintain the notation above.
The word $w=i_1^{e_1}i_2^{e_2}\cdots i_t^{e_t}\in\ooz$ is distinguished if and only if
for each $1\leq a\leq t$,
$${\widetilde {\nu(a)}}^{(i_a+1)}_{e_a}\geq {\widetilde {\nu(a)}}^{(i_a)}_1,$$
that is, the $e_a$th part of the dual partition of $\nu(a)^{(i_a+1)}$ is
greater than or equal to the first part of the dual partition of $\nu(a)^{(i_a)}$.
\end{thm}

\begin{pf} We first observe from Lemma \ref{SGE} that for each $1\leq a\leq t$
\begin{equation}\label{DISa}
M(\nu(a-1))\cong\underbrace{S_{i_a}*\cdots*S_{i_a}}_{e_a}*M(\nu(a))\cong (e_aS_{i_a})*M(\nu(a)).
\end{equation}
Thus, $M:=M(\nu(0))$ admits a filtration
\begin{equation}\label{DISb}
M=M_0\supset M_1\supset\cdots\supset M_{t-1}\supset M_t=0
\end{equation}
such that $M_{a}\cong M(\nu(a))$ and
$M_{a-1}/M_a\cong e_aS_{i_a}$ $\forall a$. This is a reduced
filtration of $M$ of type $w$.

We claim that a submodule $N_a$ of
$M(\nu(a-1))$ with
$N_a\cong M(\nu(a))$
 is unique if and only if ${\widetilde {\nu(a)}}^{(i_a+1)}_{e_a}\geq {\widetilde {\nu(a)}}^{(i_a)}_1.$
Indeed, from the definition of $\sz_{i_a}^+$, we see that
${\widetilde {\nu(a-1)}}^{(i_a)}$ is obtained by adding parts
${\widetilde {\nu(a)}}^{(i_a+1)}_j+1$ ($1\leq j\leq e_a$)
to ${\widetilde {\nu(a)}}^{(i_a)}$.
Thus, by Lemma \ref{DISL}, $N_a$ has the described property
if and only if
$$\aligned
&{\widetilde {\nu(a-1)}}^{(i_a)}_j={\widetilde {\nu(a)}}^{(i_a+1)}_j+1,
1\leq j\leq e_a,\text{ and}\cr
&{\widetilde {\nu(a-1)}}^{(i_a)}_{e_a}>{\widetilde {\nu(a-1)}}^{(i_a)}_{e_a+1}
={\widetilde {\nu(a)}}^{(i_a)}_1,
\endaligned$$
which is equivalent to the condition
${\widetilde {\nu(a)}}^{(i_a+1)}_{e_a}\geq {\widetilde {\nu(a)}}^{(i_a)}_1.$

We now prove the theorem by induction on $t$. If $t=0$ or $1$, the
theorem is obviously true. Let now $t>1$.

We first assume that $w$ is distinguished. Then the above filtration
(\ref{DISb}) is the unique
reduced filtration of $M$ of type $w$. Thus $M_1$ is the unique submodule
of $M$ such that $M_1\cong M(\nu(1))$ and $M/M_1\cong e_1S_{i_1}$. This
implies by the claim that
${\widetilde {\nu(1)}}^{(i_1+1)}_{e_1}\geq {\widetilde {\nu(1)}}^{(i_1)}_1.$
Certainly, the subword $w_1$ is also distinguished. By induction,
we infer that
${\widetilde {\nu(a)}}^{(i_a+1)}_{e_a}\geq {\widetilde {\nu(a)}}^{(i_a)}_1\;\;
\text{for $1\leq a\leq t$}.$

Conversely, assume that ${\widetilde {\nu(a)}}^{(i_a+1)}_{e_a}\geq {\widetilde {\nu(a)}}^{(i_a)}_1$
for all $1\leq a\leq t$. We also assume that $k$ is an algebraically closed
field. Let
$$M=M'_0\supset M'_1\supset\cdots\supset M'_{t-1}\supset M'_t=0$$
be any reduced filtration of $M$ of type $w$.
Then, each $M'_i$ ($0\leq i\leq t-1$) is an extension of $e_{i+1}S_{i+1}$
by $M'_{i+1}$. Inductively, we may assume $M'_{i+1}\led M_{i+1}$.
Thus, (\ref{DISa}) and Lemma \ref{GE} imply $M'_i\led M_i$.
In particular, we have $M'_1\led M_1$.
On the other hand, since $M/M'_1\cong e_1S_{i_1}$, Lemma \ref{DISL} implies
$M_1\led M'_1$. Hence, $M'_1\cong M_1\cong  M(\nu(1))$. This implies
$M_1'=M_1$ since by the claim $M$ has a unique submodule isomorphic to
$M(\nu(1))$. Now, by induction, $w_1$ is distinguished, that is,
$M_1\cong M(\nu(1))$ has a unique reduced filtration of type $w_1$.
Therefore, $M$ has a unique reduced filtration of type $w$,
i.e., $w$ is distinguished.
\end{pf}

\section{A comparison of order relations}

For $\lz,\mu\in\Pi$,
we write $\mu\led \lz$ if $M(\mu)\led M(\lz)$. Then
$\led$ defines a partial order on the set $\Pi$ which is independent
of $k$ by Lemma \ref{HOM} (and the remark afterwards).
We first
prove in this section that this order coincides with the extension order
$\peq$ introduced in \S2.

\begin{prop} \label{order} The order relations $\led$ and $\peq$ on $\Pi$ coincide.
\end{prop}

\begin{pf} For $\lz, \mu\in\Pi$, it is
easy to see from Lemma \ref{BOD} that $\mu\peq\lz$ implies $\mu\led\lz$.
It remains to prove that $\mu\led\lz$ implies $\mu\peq\lz$.
We assume that $k$ is algebraically closed and write $M(\mu)\peq M(\lz)$
if $\mu\peq\lz$.

Let $\mu\led\lz$ in $\Pi$ and set $N:=M(\mu), M:=M(\lz)$.
By Lemma \ref{HOM}(4), to prove $N\peq M$, it suffices to show that if there
is an exact sequence
\begin{equation}\label{pbkc}
0\lra X\lra M\lra Y\lra 0,
\end{equation}
then $X\oplus Y\peq M$. We apply induction on $d:=\dim_k M+\dim_k Y$.
If $d=0$ or $1$, this is clearly true. Suppose $d>1$. If $Y$ is decomposable,
say $Y=Y'\oplus Y''$ with $Y'\not=0\not=Y''$, then we obtain a split exact
sequence $0\ra Y'\ra Y\ra Y''\ra 0$. This together with (\ref{pbkc})
gives a commutative diagram as described in Lemma \ref{pbk}. In particular,
we have exact sequences
$$
0\ra X\ra M'\ra Y'\ra 0\text{ and }
0\ra M'\ra M\ra Y''\ra 0
$$
Since both $\dim_k M'+\dim_k Y'<d$ and $\dim_k M+\dim_k Y''<d$, the inductive
hypothesis implies
$$X\oplus Y'\peq M'\;\;\text{and}\;\;M'\oplus Y''\peq M.$$
This concludes $X\oplus Y=X\oplus Y'\oplus Y''\peq M'\oplus Y''\peq M$. Hence,
we may suppose that $Y$ is indecomposable. Similarly, $X$ can also be supposed
to be indecomposable. In the case where both $X$ and $Y$ are indecomposable,
$X\oplus Y\peq M$ follows from Lemma \ref{INEXT}.
\end{pf}

With the coincidence of the two orders, results 5.4 and Theorem C
in \cite{R93}
continue to hold if the order $\peq$ is replaced by $\led$.
However, if we use the order $\led$ at the outset, we may use the
nice properties of generic extensions.
The  following ``dominance'' property is given in  \cite[5.4]{R93}.
Here we provide a new and short proof.

\begin{lem} \label{DOM}Let $w\in\ooz$ and $\lz\ged\mu$ in $\Pi$. Then
$\lr{w, M(\lz)}\not=0$ implies $\lr{w,M(\mu)}\not=0$.
\end{lem}

\begin{pf} In view of Lemma \ref{COM}, we may suppose that
$k$ is an algebraically closed field.

Let $w=i_1i_2\cdots i_m$ and set
$w'=i_2\cdots i_m$. We apply induction on $m$. If $m=1$ then $\lz\ged\mu$
forces $M(\lz)=M(\mu)$ and the result is clear. Let now $m>1$.
Since $\lr{w,M(\lz)}\not=0$, $M(\lz)$ has a submodule $M'$ such that
$\lr{w',M'}\not=0$ and $M(\lz)/M'\cong S_{i_1}$. Thus, by Prop. \ref{GE},
we obtain $M(\lz)\led S_{i_1}\ast M'.$ This together with $M(\mu)\led M(\lz)$
yields $M(\mu)\led S_{i_1}\ast M'.$
Now, applying Prop. \ref{GE} to this relation,
there are modules $N', N''$ with $N'\led M'$,
$N''\led S_{i_1}$, and an exact sequence
$$0\lra N'\lra M(\mu)\lra N''\lra 0.$$
Since $S_{i_1}$ is simple, $ N''\led S_{i_1}$ implies $N''\cong S_{i_1}$.
On the other hand, the inductive hypothesis implies
$\lr{w',N'}\neq0$, that is, $N'$ has a composition series
of type $w'$. Therefore, $M(\mu)$ has a composition series of type $w$,
i.e., $\lr{w,M(\mu)}\not=0$.
\end{pf}

Theorem C in \cite{R93} plays a key role in the proofs of the main
results there. It holds only for condensed words and requires a rather
long proof. We are now able to generalize this result by
removing the restriction on the words.

\begin{thm} \label{DOMI}For each $w\in\ooz$ and each $\pi\in\Pi$,
$\lr{w, M(\pi)}\not=0$ if and only if $\pi\led\wp(w)$.
\end{thm}

\begin{pf} Again by Lemma \ref{COM} and Remark \ref{GEk}(2), we may assume that $k$ is an
algebraically closed field.

Let $w=i_1i_2\cdots i_m\in\ooz$, $\pi\in\Pi$, and
$\lr{w, M(\pi)}\not=0$. We use induction on $m$ to prove that $\pi\led\wp(w)$,
i.e., $M(\pi)\leq M(\wp(w))$. If $m=0$, or $1$, there is nothing to prove.
Let $m>1$ and set $w'=i_2\cdots i_m$. Then
$$M(w)=S_{i_1}\ast (S_{i_2}\ast\cdots\ast S_{i_m})=S_{i_1}\ast M(w').$$
Since $\lr{w,M(\pi)}\not=0$, $M(\pi)$ has a submodule $M'$ with
$M(\pi)/M'\cong S_{i_1}$ and $M'$ has a composition series of type $w'$.
Thus, there is an exact sequence
$$0\lra M'\lra M(\pi)\lra S_{i_1}\lra 0.$$
By the inductive hypothesis, we have $ M'\led M(w')$. We infer from
Prop. \ref{GE} that
$$M(\pi)\led S_{i_1}\ast M'\led S_{i_1}\ast M(w')=M(w)=M(\wp(w)),$$
that is, $ \pi\led \wp(w)$.

Conversely, let $\pi\led \wp(w)$. By the definition of $\wp(w)$,
$\lr{w, M(\wp(w))}\not=0$. This together with Lemma \ref{DOM} implies
$\lr{w,M(\pi)}\not=0$.
\end{pf}

For $\lz\in\Pi$, let $\Pi^{\leq\lz}$ be the poset ideal generated by
$\lz$, i.e.,
$$\Pi^{\leq\lz}=\{\mu\in\Pi\mid \mu\leq\lz\}.$$

\begin{cor} If $\lz\in\Pi^s$, then $\mu\in\Pi^{\leq\lz}$ if and only if
there exists a $w\in\wp^{-1}(\lz)$ such that $M(\mu)$ has a composition series
of type $w$.
\end{cor}
It would be interesting to describe the ideal $\Pi^{\leq\lz}$
for an arbitrary $\lz\in\Pi$.

\section{Hall algebras and its composition subalgebra}

Given three modules $L, M, N$ in $\bbt$, let $F^L_{MN}$ be the
number of submodules $V$ of $L$ such that $V\cong N$ and $L/V\cong
M$. More generally, given modules $M, N_1,\cdots,N_m$ in
$\bbt$, we let $F^M_{N_1\,\cdots\,N_m}$ be the number of the
filtrations
$$M=M_0\supseteq M_1\supseteq \cdots \supseteq M_{m-1}\supseteq M_m=0,$$
such that $M_{t-1}/M_t\cong N_t$ for all $1\le t \le m$.
In particular, we have, for $w=i_1\cdots i_m\in\ooz$,
$F^M_{S_{i_1}\cdots S_{i_m}}=\lr{w,M}.$
Note that
\begin{equation}\label{assoc}
F^M_{N_1N_2N_3}=\sum_{[N]}F_{N_1N}^M F_{N_2N_3}^N.
\end{equation}

By \cite{R93} and \cite{Guo}, there exist Hall polynomials in $\bbt$. More precisely,
for $\pi,\mu_1,\cdots,\mu_m$ in $\Pi$, there is a polynomial
$\vphi^\pi_{\mu_1\cdots \mu_m}(q)\in\ca:=\bbz[q]$
such that for any finite $k$ of $q_k$ elements
$$\vphi^\pi_{\mu_1\cdots\mu_m}(q_k)=F^{M_k(\pi)}_{M_k(\mu_1)\cdots M_k(\mu_m)}.$$
In particular, when $M(\mu_j)=S_{i_j}$ for all $j$, we denote the polynomial
$\vphi^\pi_{\mu_1\cdots\mu_m}$ by $\lr{w|\pi}$ as in \cite[8.1]{R93}.

\begin{lem}\label{DISC} Let $w\in\ooz$ be a distinguished word with the
tight form $i_1^{e_1}i_2^{e_2}\cdots i_t^{e_t}$. Then
$$\lan w|\wp(w)\ran=\prod_{a=1}^t[\![e_a]\!]^!,$$
where $[\![e_a]\!]^!=[\![1]\!]\cdots [\![e_a]\!]$ with
$[\![m]\!]=\frac{1-q^m}{1-q}\in\ca.$
\end{lem}

\begin{pf} We first observe that, since every composition series
of type $w$ is a refinement of the associated reduced filtration of
type $w$,
$$\lr{w,M(w)}=F_{e_1S_{i_1}\cdots e_tS_{i_t}}^{M(w)}
\prod_{a=1}^t\lr{i_a^{e_a},{e_aS_{i_a}}}.$$
If $w$ is distinguished, it follows that $M(w)$ has
a unique reduced filtration for an infinite number of finite
fields. Thus $F_{e_1S_{i_1}\cdots e_tS_{i_t}}^{M(w)}\equiv 1$
(cf. Lemma \ref{COM}). Now the result follows from
the fact that the number of composition series
of $e_aS_i$ is $[\![e_a]\!]^!$ (cf. \cite[8.2]{R93}).
\end{pf}

 The generic (untwisted)
{\it Ringel-Hall algebra} $\ch_q(\dt)$
of $\dt$ is by definition the free $\ca$-module with basis $\{u_\pi|\pi\in\Pi\}$,
and the multiplication is given by
$$u_\mu u_\nu=\sum_{\pi\in\Pi}\vphi^\pi_{\mu\nu}(q)u_\pi.$$
In practice, we sometimes write $u_\pi=u_{[M(\pi)]}$ in order to make certain
calculations in term of modules. For $i\in I$, we set $u_i=u_{[S_i]}$
and denote by $\cc_q(\dt)$ the subalgebra of $\ch_q(\dt)$ generated by $u_i$,
$i\in I$. This is
 called the (generic) {\it composition algebra} of $\dt$. It is easy
to see that $\cc_q(\dt)$ is a proper subalgebra of $\ch_q(\dt)$.
Moreover, both $\ch_q(\dt)$ and $\cc_q(\dt)$ admit a natural $\bbn^n$-grading
by dimension vectors.

Let $\ca\ooz$ be the semigroup algebra of $\ooz$, which is indeed the free
$\ca$-algebra generated by the set $I$. Thus, we have an algebra homomorphism
$f:\ca\ooz\to\cc_q(\dt)$ with $f(i)=u_i$ for all $i\in I$.

Following \cite[1.2]{R93},
there is an $\ca$-bilinear form
$$\lan-|-\ran:\ca\ooz\times \ch_q(\dt)\ra\ca,\;(w,u_\pi)\lmto\lan w|\pi\ran(q).$$
We denote by $\fr$ 
the set of $x\in\ca\ooz$ 
such that $\lan x|-\ran=0$. 
We call $\fr$ the {\it left radical} of the form.

For each word $w=i_1i_2\cdots i_m\in\ooz$, we define
$$u_w=u_{i_1}u_{i_2}\cdots u_{i_m}\in\cc_q(\dt).$$

\begin{lem} \label{CA}
We have $\fr=\ker(f)$. Thus $f$ induces an isomorphism
$\cc_q(\dt)\cong\ca\ooz/\fr$. Moreover,
$\cc_q(\dt)\otimes_{\bbz[q]}\bbq[q]$ is free.
\end{lem}

\begin{pf} Define the $\ca$-bilinear form
$$(-,-):\ch_q(\dt)\times \ch_q(\dt)\ra\ca$$
by $(u_\lz,u_\mu)=\delta_{\lz\mu}$. Then
$\lr{w|\pi}=(u_w,u_\pi)$. Thus, for any $x=\sum_w x_ww\in\ca\ooz$,
we have
$$\lr{x|\pi}=\sum_wx_w(u_w,u_\pi)=(f(x),u_\pi).$$
Therefore, the equality follows.\end{pf}

For any commutative ring $\ca'$ which is an $\ca$-algebra and
any $\ca$-module $M$, let $M_{\ca'}= M\otimes_\ca\ca'$
denote the $\ca'$-module obtained from $M$
by base change to $\ca'$. Clearly, base change
induces an $\ca'$-bilinear form
$$\lan-|-\ran_{\ca'}:\ca'\ooz\times \ch_q(\dt)_{\ca'}\ra\ca'.$$
It is easy to see that, if $\ca'$ is an integral domain containing
$\ca$, then $\fr_{\ca'}$ is the same as the left radical of
$\lan-|-\ran_{\ca'}$.

Ringel discovered that the ideal $\fr$ contains the ideal
$\fr'$ generated
by those near quantum Serre relations (see \cite[8.6]{R93}) and that
these two ideals are equal after base change
to the localization $\ca_{(q-1)}=\bbz[q]_{(q-1)}=\bbq[q]_{(q-1)}$
of $\ca$ at the maximal ideal generated by $q-1$.
Thus, we have the following theorem
(see \cite[8.7]{R93}).

\begin{thm}\label{CREL} Let $\ca'=\ca_{(q-1)}$.

(1) If $n\geq 3$, then $\cc_q(\dt_n)_{\ca'}$ is generated by $u_i$, $i\in I$,
with relations:
$$\aligned
&u_i u_j=u_j u_i\;\;\mbox{if $j\not\equiv i\pm 1$ (mod $n$),
$i,j\in I$},\cr
&u_i^2 u_{i+1}-(q+1) u_i u_{i+1}u_i +q u_{i+1} u_i^2=0,\cr
&u_i u_{i+1}^2-(q+1) u_{i+1}u_i u_{i+1}+q u_{i+1}^2u_i=0,\;
i\in I.
\endaligned$$

(2) The algebra $\cc_q(\dt_2)_{\ca'}$ is generated by $u_1$ and $u_2$
with relations:
$$\aligned
&q u_1^3u_2-(q^2+q+1) u_1^2 u_2 u_1+(q^2+q+1)u_1 u_2 u_1^2-q u_2 u_1^3=0,\cr
&q u_2^3u_1-(q^2+q+1) u_2^2 u_1 u_2+(q^2+q+1)u_2 u_1 u_2^1-q u_1 u_2^3=0.
\endaligned$$
\end{thm}

The following results generalize \cite[8.5, 8.8]{R93}.

\begin{thm} \label{MBC} For every $\pi\in\Pi^s$, choose an arbitrary word
$w_\pi\in\wp^{-1}(\pi)$.

(1) The free $\ca$-submodule $\frak C$ of $\ca\ooz$ spanned by
$\{w_\pi|\pi\in\Pi^s\}$ intersects $\fr$ trivially.

(2)  We have $\bbq(q)\ooz={\frak C}_{\bbq(q)}\oplus\fr_{\bbq(q)}.$

(3)  Let $\ca'=\ca_{(q-1)}$ and assume that all $w_\pi$
are distinguished. Then
$\ca'\ooz={\frak C}_{\ca'}\oplus\fr_{\ca'}.$

\end{thm}

\begin{pf} If $x=\sum_{i=1}^ma_iw_i\in{\frak C}\cap \fr$, then
$\wp(w_i)\neq \wp(w_j)$ for $i\neq j$. It follows from
 Theorem \ref{DOMI} that
the matrix $(\lr{w_i|\wp(w_j)})$ is upper triangular
under an appropriate ordering and the diagonal entries are non-zero.
This implies $x=0$,
proving (1).

To prove (2) resp. (3), it remains
to prove that $\bbq(q)\ooz={\frak C}_{\bbq(q)}+\fr_{\bbq(q)}$
resp. $\ca'\ooz={\frak C}_{\ca'}+\fr_{\ca'}.$
This can be done  in a way similar to the proof of
\cite[8.8]{R93}, using
Theorems \ref{DOMI} and \ref{CREL}.
Note that, if $w_\pi$ is distinguished, then, by Lemma \ref{DISC},
$\lan w|\wp(w)\ran$ is invertible in $\ca'$. This is just the case
considered by Ringel.
\end{pf}

An immediate consequence is the following monomial basis theorem for
the composition algebra $\cc_q(\dt)_{\ca'}$.

\begin{cor}\label{MON}  For every $\pi\in\Pi^s$, choose an arbitrary word
$w_\pi\in\wp^{-1}(\pi)$. The set $\{u_{w_\pi}|\pi\in\Pi^s\}$ is a $\bbq(q)$-basis
of $\cc_q(\dt)_{\bbq(q)}$. Moreover, if all $w_\pi$ are
chosen to be distinguished, then this set is an $\ca'$-basis of
$\cc_q(\dt)_{\ca'}$ where $\ca'=\ca_{(q-1)}$.
\end{cor}


\section{Proof of Theorem 1.1}

We now in this section transfer
the monomial bases for the composition algebra of a cyclic quiver
 to a quantum affine ${\frak sl}_n$.
We need to consider the twisted version of Hall and composition algebras
(see \cite{R932}).

First, we recall that the Euler form associated with the cyclic quiver $\dt$
is the bilinear form
$\vez(-,-):\bbz^n\times\bbz^n\ra \bbz$ defined by
$$\vez({\bf a}, {\bf b})=\sum_{i=1}^n a_i b_i-
\sum_{i=1}^{n} a_i b_{i+1},$$
where ${\bf a}=(a_1, \cdots, a_n)$, ${\bf b}=(b_1, \cdots, b_n)$, and
$b_{n+1}=b_1$. It is well-known that for two representations $M,N\in\bbt$,
there holds
$$\vez(\udim M, \udim N)=\dim_k\hom(M, N)-\dim_k\ext^1(M, N).$$

Let $\cz=\bbz[v,v^{-1}]$, where $v$ is an indeterminate with $v^2=q$.
The {\it twisted Ringel-Hall algebra} $\ch_v^\star(\dt)$ of $\dt$ is
by definition the
free $\cz$-module with basis $\{u_\pi=u_{[M(\pi)]}|\pi\in\Pi\}$, and
the multiplication is defined by
$$u_\mu\star u_\nu=v^{\vez(\mu,\nu)}\sum_{\pi\in\Pi}
\vphi^{\pi}_{\mu\nu}(v^2)u_\pi,$$
where $\vez(\mu,\nu)=\vez(\udim M(\mu),\udim M(\nu))$. The $\cz$-subalgebra
$\cc_v^\star(\dt)$
of $\ch_v^\star(\dt)$ generated by $u_i$, $i\in I$, is called
the {\it twisted composition algebra}.

Let $U^+$ be the $\cz$-subalgebra
of $\bU=\bU_v(\widehat {\frak sl}_n)$ generated by $E_i,i\in I$,
and let $\cz'=\cz_{(v-1)}$ be the localization of $\bbz[v,v^{-1}]$
at the ideal generated by $v-1$. Since the relations in
Theorem \ref{CREL} become the quantum Serre relations under the twisted
multiplication, we obtain by modifying Ringel's proof of
\cite[8.7]{R93} an isomorphism
$$\Phi:\cc_v^\star(\dt)_{\cz'}\cong U^+_{\cz'},\;u_i\mapsto E_i,\;i\in I.$$

For each $w=i_1i_2\cdots i_m\in\ooz$, if we put
$$E_w=E_{i_1}E_{i_2}\cdots E_{i_m}\in U^+$$
as in the introduction, we have the following (cf. Theorem
\ref{MON}).

\begin{thm} \label{MBTa} For every $\pi\in\Pi^s$, choose an arbitrary word
$w_\pi\in\wp^{-1}(\pi)$. The set $\{E_{w_\pi}|\pi\in\Pi^s\}$ is a $\bbq(v)$-basis
of $\bU^+$. Moreover, if all $w_\pi$ are
chosen to be distinguished, then this set is a $\cz'$-basis of
$U^+_{\cz'}$.
\end{thm}
This theorem together with the triangular decomposition of $\bU$
gives Theorem \ref{MBT}.

Let ${\frak g}={\frak n}_-\oplus {\frak h}\oplus {\frak n}_+$ be
the affine Lie algebra $\widehat{\frak sl}_n$ over $\bbq$
of type $\tilde A_{n-1}$ with generators $e_i, f_i, h_j$.
Let ${\frak U}(\frak g)$ be the universal enveloping algebra of ${\frak g}$.
We also define monomials $e_w$ similarly for $w\in\ooz$ in
${\frak U}({\frak n}_+)$. Then, we have the following.

\begin{cor}  For every $\pi\in\Pi^s$, choose an arbitrary
distinguished word
$w_\pi\in\wp^{-1}(\pi)$. The set $\{e_{w_\pi}|\pi\in\Pi^s\}$ is a $\bbq$-basis
of ${\frak U}({\frak n}_+)$.
\end{cor}
\begin{pf}
By Theorem \ref{CREL}, we have $\cc_q(\dt)_{\ca'}/(q-1)\cc_q(\dt)_{\ca'}\cong
{\frak U}({\frak n}_+)$. The result follows from Theorem
\ref{MON}.
\end{pf}

\section{A comparison of bases}

With the isomorphism $\Phi$ above, we now identify $\bU^+$ with
$\cc_v^\star(\dt)_{\bbq(v)}$. Thus, $E_i=u_i$ and, for
each $w=i_1i_2\cdots i_m\in\ooz$, we have
$$E_w=u_w^\star:=u_{i_1}\star u_{i_2}\star\cdots\star u_{i_m}
=v^{\vez(w)}u_w,$$
where $\vez(w)=\sum_{r<t}\vez(\udim S_{i_r},\udim S_{i_t})$.
Further, we have the following relation between monomial bases
$\{E_{w_\pi}\}_{\pi\in\Pi^s}$
and the defining basis $\{u_\lz\}_{\lz\in\Pi}$ of $\ch_v^\star(\dt)$.

\begin{prop} \label{UW} For each $w\in\ooz$, we have
\begin{equation}\label{UWa}
E_w=\sum_{\lz\led\wp(w)}v^{\vez(w)}\lan w|\lz\ran u_\lz.
\end{equation}
Moreover, the coefficients appearing in the sum are all non-zero and all
in $\bbn[v,v^{-1}]$.
\end{prop}

\begin{pf} By Theorem \ref{DOMI}, it remains to show that each
polynomial $\lan w|\lz\ran\in\bbz[q]$ has non-negative coefficients.

Let $w=i_1i_2\cdots i_m$. If $m=1$ or $2$, $\lan w|\lz\ran\in\bbz[q]$ has
obviously the required property. Assume $m\geq 2$ and set $w'=i_2\cdots i_m$.
Then we have by (\ref{assoc})
$$\lan w|\lz\ran=\sum_\mu \vphi^\lz_{\nu\mu}(q)\lan w'|\mu\ran,$$
where $\nu\in\Pi$ satisfies $M(\nu)=S_{i_1}$. By induction, it suffices
to prove that each $\vphi^\lz_{\nu\mu}(q)$ lies in $\bbn[q]$.
Let $k$ be a finite field of $q_k$ elements. We now
calculate the number $f:=F^{M_k(\lz)}_{M_k(\nu)M_k(\mu)}$.
Since $M_k(\nu)=S_{i_1}$, we have clearly $f=F^M_{S_{i_1}N}$ with
$$M=M_{i_1}(\tlz^{(i_1)})=\oplus_{j}S_{i_1}[\tlz^{(i_1)}_j] \text{ and }
N=M_{i_1}(\tilde\mu^{(i_1)}).$$
For simplicity, we write $i=i_1$, $p:=\tlz^{(i_1)}=(p_1,\cdots, p_l)$.
By Lemma \ref{EXT}(2), 
$F^M_{S_iN}\neq0$ implies $N\cong N_t$
for some $1\leq t\leq l$, where
$N_t=S_{i+1}[p_t-1]\oplus\bigoplus_{j:j\not=t} S_i[p_j]$.
Let $a$ (resp. $b$, $c$) be the number
of $p_j$'s such that $p_j>p_t$ (resp. $p_j=p_t$, $p_t>p_j$) so that
$l=a+b+c$. We claim
that
$$f=F^M_{S_iN_t}=q_k^a(1+q_k+\cdots q_k^{b-1})=q_k^a[\![b]\!].$$
Indeed, choose $v_j\in S_i[p_j]\backslash \rad(S_i[p_j])$ for each
$1\leq j\leq l$. Then (the images of) $v_1,\cdots,v_l$ form a basis
for $M/\rad(M)$. Let $V$  denote the subspace of $M$ spanned by
$v_{a+1},\cdots,v_{a+b}$ and let $\Xi$ be
the set of subspaces of $V$ of codimension 1. Then $|\Xi|=[\![b]\!]$.
For each $W\in\Xi$ and $(x_1,\cdots x_a)\in k^a$, we fix
an element $v_W\in V\backslash W$ and define $N(W, x_1,\cdots x_a)$ to be
the submodule of $M$ generated by $v_1+x_1v_W,\cdots, v_a+x_av_W,
W, v_{a+b+1},\cdots v_l$ and $\rad(M)$.
Clearly, the modules $N(W, x_1,\cdots x_a)$ are all distinct
submodules isomorphic to $N_t$ and
each submodule $N$ of $M$ with $N\cong N_t$ is of this form.  Therefore,
$f=|\Xi\times k^a|=q_k^a[\![b]\!],$ as required.
Consequently, each $\vphi^\lz_{\nu\mu}(q)$ lies in $\bbn[q]$.
\end{pf}

We observe that the relation (\ref{UWa}) shares similar properties
possessed by the canonical basis  $\{b_\pi\}_{\pi\in\Pi}$
introduced in \cite{L91}.
It is known (cf. \cite[(8)]{LTV}) that, when writing $b_\pi$ as
a linear combination of the defining basis $\{u_\lz\}_{\lz\in\Pi}$,
both order and positivity properties are satisfied.
More precisely, we have
\begin{equation}\label{UWb}
b_\pi=\sum_{\lz\leq\pi}p_{\lz,\pi}u_\lz, \,\,\,p_{\lz,\pi}\in\bbn[v,v^{-1}].
\end{equation}
However, there is no elementary proof for the positivity property in this case.
 Note that the subset
$\{b_\pi\}_{\pi\in\Pi^s}$ forms a basis for $\bU^+$.
It is also observed (see \cite[p.40]{LTV}) that non-separated $\pi$ may
occur in the right hand side of (\ref{UWb}) (and of (\ref{UWa}) as well) .


Motivated by \cite[Thm 4]{R93}, we now construct another basis for $\bU^+$
from the defining basis $\{u_\lz\}_{\lz\in\Pi}$.
Recall from \S 7 the $\bbq(v)$-bilinear form
$\lan-|-\ran:\bbq(v)\ooz\times \ch_q(\dt)_{\bbq(v)}\ra\bbq(v)$.
Let $\fs$ denote the right radical of this form, that is, the space
of all $y\in\ch_q(\dt)_{\bbq(v)}$ with $\lr{-|y}=0$. Obviously,
$\fs=\oplus_{\bf d}\fs_{\bf d}$ admits
the grading by dimension vectors. By \cite[Thm 4]{R93}, we
see that $\fs$ is a direct complement of the subspace of
$\ch_q(\dt)_{\bbq(v)}$ spanned by all $u_\pi,\pi\in\Pi^s$.
However, in contrast to the left radical $\fr$, $\fs$ is in general not
an ideal of $\ch_q(\dt)_{\bbq(v)}$. For example, let $n=2$. Then
$\fs_{(1,1)}$ is one-dimensional over $\bbq(v)$
with a basis element $x=-u_\lz+u_\mu+u_\nu$, where
$\lz=\bigl((1),(1)\bigr)$, $\mu=\bigl((1,1),\emptyset\bigr)$, and
$\nu=\bigl(\emptyset,(1,1)\bigr)$ are pairs of
partitions in $\Pi$. But an easy calculation shows that neither $u_1x$ nor $xu_1$
lies in $\fs_{(2,1)}$. Thus $\fs$ is not an ideal of $\ch_q(\dt_2)_{\bbq(v)}$.

In order to obtain an ideal, we use Green's form \cite{Gr} which is a
modified version of the bilinear form
$\lr{-|-}$. We consider the twisted version $\ch_v^\star(\dt)_{\bbq(v)}$ and
$\cc_v^\star(\dt)_{\bbq(v)}=\bU^+$.

For each $\pi\in\Pi$, it is well-known that there is a monic polynomial
$a_\pi\in\bbz[q]$ such that $a_\pi(q_k)=|\mbox{Aut}(M_k(\pi))|$ for any finite
field $k$ of $q_k$ elements.
Define a $\bbq(v)$-bilinear from
$$\lr{-|-}':\bbq(v)\ooz\times\ch_v^\star(\dt)_{\bbq(v)}\lra\bbq(v)$$
by setting $\lr{w|u_\pi}'=v^{\vez(w)}a_\pi^{-1}\lr{w|u_\pi}$, and denote
by $\fs'$ the right radical subspace of this form.
Let $\bbq(v)\Pi^s$ be the subspace of $\ch_v^\star(\dt)_{\bbq(v)}$ spanned
by all $u_\pi,\pi\in\Pi^s$.

\begin{lem}\label{IDEAL} The subspace $\fs'$ is an ideal of $\ch_v^\star(\dt)_{\bbq(v)}$
and moreover
$$\ch_v^\star(\dt)_{\bbq(v)}=\bbq(v)\Pi^s\oplus\fs'.$$
\end{lem}

\begin{pf} Following \cite{Gr}, there is
a symmetric and non-degenerate bilinear form
$$(-,-)':\ch_v^\star(\dt)_{\bbq(v)}\times \ch_v^\star(\dt)_{\bbq(v)}\lra\bbq(v)$$
defined by $(u_\lz,u_\mu)'=a_\lz^{-1}\dz_{\lz\mu}$. Then
$\lr{w|u_\pi}'=(u_w^\star,u_\pi)'$. This implies
$$\fs'_\bd=\{y\in {(\ch_v^\star(\dt)_{\bbq(v)})}_\bd|
({(\cc_v^\star(\dt)_{\bbq(v)})}_\bd, y)'=0\}$$
for each $\bd\in\bbn^n$.
By \cite[Thm 1]{Gr}, we see easily that $\fs'$ is an ideal of
$\ch_v^\star(\dt)_{\bbq(v)}$.

Now, by an argument similar to the proof of \cite[Thm 4]{R93},
we see that $\bbq(v)\Pi^s\cap \fs'=0$. On the other hand,
since $(-,-)'$ is non-degenerate and its restriction to
$\cc_v^\star(\dt)_{\bbq(v)}$ is again non-degenerate
(see \cite{Gr} and \cite[1.2, 33.1]{L93}), we obtain
$(\ch_v^\star(\dt)_{\bbq(v)})_\bd=(\cc_v^\star(\dt)_{\bbq(v)})_\bd
\oplus \fs'_\bd$.  Thus,
$$\dim\fs'_\bd=\dim (\ch_v^\star(\dt)_{\bbq(v)})_\bd
-\dim (\cc_v^\star(\dt)_{\bbq(v)})_\bd=|\Pi_\bd|-|\Pi^s_\bd|,$$
where $\Pi_\bd=\{\pi\in\Pi|\udim M(\pi)=\bd\}$ and
$\Pi^s_\bd=\Pi^s\cap\Pi_\bd$. Hence,
$$\ch_v^\star(\dt)_{\bbq(v)}=\bbq(v)\Pi^s\oplus\fs'.$$
\end{pf}

\begin{thm}\label{IDEALC} We have $\bU^+\cong \ch_v^\star(\dt)_{\bbq(v)}/\fs'$.
Moreover, the algebra $\bU^+$ has a basis
of the form $\cp:=\{\bar u_\pi\mid \pi\in\Pi^s\}$
where $\bar u_\pi$ is the homomorphic image of $u_\pi$.
\end{thm}

\begin{pf}
By the decomposition $\ch_v^\star(\dt)_{\bbq(v)}=\cc_v^\star(\dt)_{\bbq(v)}\oplus \fs'$, the required isomorphism is given by the composition
of the canonical inclusion $\cc_v^\star(\dt)_{\bbq(v)}\hookrightarrow \ch_v^\star(\dt)_{\bbq(v)}$
and the projection
$\ch_v^\star(\dt)_{\bbq(v)}
\twoheadrightarrow \ch_v^\star(\dt)_{\bbq(v)}/\fs'.$
The rest of the proof follows from the previous lemma.
\end{pf}

\begin{rem}
The basis $\cp$ given in Theorem \ref{IDEALC} may be viewed as a
PBW type basis. We remark that the transition matrix from
any monomial basis (or the canonical
basis) to the basis $\cp$ is no longer in triangular form.
That is, if we rewrite the right hand side of (\ref{UWa}) or
(\ref{UWb})
by replacing those $u_\pi$ with $\pi\not\in\Pi^s$ by
linear combinations of elements of $\cp$, then the
order relation is no longer sustained.
\end{rem}

\begin{example} Consider the case $n=2$. Then $\ooz_{(2,1)}$
consists of three words $w_1=1^22, w_2=21^2, w_3=121$, and $\Pi_{(2,1)}$
contains four elements
$$\aligned
&\az=\bigl((2),(1)\bigr),\;\bz=\bigl((2,1),\emptyset\bigr),\cr
&\gz=\bigl((1),(1,1)\bigr),\dz=\bigl((1,1,1),\emptyset\bigr),
\endaligned$$
which are ordered by $\az<\bz, \az<\gz, \bz<\dz$ and
$\gz<\dz$.
Clearly, $\wp(w_1)=\bz,\wp(w_2)=\gz$ and $\wp(w_3)=\dz$ are
separated. Further, $\fs'_{(2,1)}$ is one-dimensional with a basis element
$-(v^4-1)u_{\az}+u_{\bz}+u_{\gz}+u_{\dz}$,
that is, in $\bU^+_{(2,1)}$,
$${\bar u}_{\az}=\frac{1}{v^4-1}({\bar u}_{\bz}
+{\bar u}_{\gz}+{\bar u}_{\dz}),$$
and $\bU^+_{(2,1)}$ has bases
$\{E_{w_1},E_{w_2},
E_{w_3}\}$ and
$\{{\bar u}_{\bz},{\bar u}_{\gz},{\bar u}_{\dz}\}$
satisfying
$$\begin{array}{rl}
&E_{w_1}=\frac{1}{v(v^2-1)}
(v^4{\bar u}_{\bz}+{\bar u}_{\gz}+{\bar u}_{\dz}),\\
&E_{w_2}=\frac{1}{v(v^2-1)}
({\bar u}_{\bz}+v^4{\bar u}_{\gz}+{\bar u}_{\dz}),\\
&E_{w_3}=\frac{v}{v^2-1}
({\bar u}_{\bz}+{\bar u}_{\gz}+{\bar u}_{\dz}).
\end{array}$$
\end{example}


\begin{thebibliography}{9}

\bibitem{Bo} K. Bongartz,
{\em On degenerations and extensions of finite dimensional modules},
Adv. Math.
{\bf 121} (1996), 245-287.

\bibitem{CX} V. Chari and N. Xi, {\em Monomial bases of quantized enveloping
algebras}, In: Recent developments in quantum affine algebras and related
topics (Raleigh, NC, 1998), Contemp. Math. {\bf 248} (1999),  69--81.

\bibitem{DP} J. Du and B. Parshall, {\em
Monomial bases for $q$-Schur algebras,} Trans. Amer. Math Soc.
{\bf 355} (2003) 1593--1620.

\bibitem{Gr} J.~A. Green,
{\em Hall algebras, hereditary algebras and quantum groups},
Invent. math. {\bf 120} (1995), 361-377.

\bibitem{Guo} J.~Y. Guo,
{\em The Hall polynomials of a cyclic serial algebra},
Comm. Algebra
{\bf 23} (1995), 743-751.

\bibitem{LTV} B. Leclerc, J.-Y. Thibon and E. Vasserot,
{\em Zelevinsky's involution at roots of unity}, J. reine
angew. Math. {\bf 513} (1999), 33-51.

\bibitem{L90} G. Lusztig, {\em Canonical bases arising
from quantized enveloping algebras,} J. Amer. Math. Soc.
{\bf 3} (1990), 447-498.

\bibitem{L91} G. Lusztig, {\em Quivers, perverse sheaves, and the quantized
enveloping algebras,} J. Amer. Math. Soc.
{\bf 4} (1991), 366-421.

\bibitem{L93} G. Lusztig, {\em Introduction to Quantum
Groups,} Birkh\"auser, Boston, 1993.

\bibitem{Re} M. Reineke,
{\em Generic extensions and multiplicative bases of quantum groups
at $q=0$},
Represent. Theory {\bf 5} (2001), 147-163.

\bibitem{Re1} M. Reineke,
{\em Feigin's map and monomial bases for quantized enveloping
algebras}, Math. Z. {\bf 237} (2001), 639-667.

\bibitem{R90} C.~M. Ringel,
{\em Hall algebras and quantum groups},
Invent. math.
{\bf 101} (1990), 583-592.

\bibitem{R93} C.~M. Ringel,
{\em The composition algebra of a cyclic quiver},
Proc. London Math. Soc.
{\bf 66} (1993), 507-537.

\bibitem{R932} C.~M. Ringel,
{\em Hall algebras revisited},
Israel Mathematical Conference Proceedings,
Vol. {\bf 7} (1993), 171-176.

\bibitem{R95} C.~M. Ringel,
{\em The Hall algebra approach to quantum groups},
Aportaciones Matem\'aticas Comunicaciones
{\bf 15} (1995), 85-114.


\bibitem{Zwa} G. Zwara,
{\em Degenerations for modules over representation-finite biserial algebras},
J. Algebra
{\bf 198} (1997), 563-581.



\end{thebibliography}
\end{document}